\documentclass[a4paper,10pt]{article}
\usepackage{amsfonts}
\usepackage{graphicx,amssymb,amsmath}
\usepackage{amscd}
\usepackage[english]{babel}
\usepackage[ansinew]{inputenc}
\usepackage{slashbox}
\usepackage[all]{xy}
\usepackage{latexsym}
\usepackage{amsmath}
\usepackage{stmaryrd}

 \hsize \textwidth
 \advance \hsize by -\marginparwidth
\newcommand{\field}[1]{\mathbb{#1}}
\newcommand{\R}{\field{R}}

\newcommand{\N}{\field{N}}

\newtheorem{theo}{Theorem}[section]
\newtheorem{theor}{Theorem}

\newtheorem{defi}[theo]{Definition}
\newtheorem{lem}[theo]{Lemma}
\newtheorem{cor}[theo]{Corollary}
\newtheorem{prop}[theo]{Proposition}
\newtheorem{rem}[theo]{Remark}

\newtheorem{nota}[theo]{Notations}
\newtheorem{com}[theo]{Comments}
\newtheorem{exe}[theo]{Example}
\newtheorem{exes}[theo]{Examples}

\def\dis{\ds}
\def\ap{\rightarrow}

\def\l{\lambda}

 \def\a{\alpha}
\def\b{\beta}
\def\g{\gamma}
\def\G{\Gamma}
\def\D{\triangle}
\def\t{\tau}
\def\d{\delta}
\def\o{\omega}

\def\l{\lambda}
\def\L{\Lambda}
\def\n{\nu}

\def\r{\rho}

\def\m{\mu}
\def\n{\nu}
\def\s{\sigma}

\def\so{\underline}

\def\O{\Omega}

\def\e{\epsilon}

\def\~{\tilde}
\def\dis{\displaystyle}
\input epsf

\title{ On Sussmann theorem for orbits of  sets of  vector fields on  Banach manifolds  }
\author{  A. Lathuille \& F. Pelletier\footnote
{Laboratoire de Math\'ematiques, Universit\'e de Savoie, Campus scientifique, 73376 Le Bourget du Lac Cedex} }
\date{}

\begin{document}
\maketitle

\begin{abstract}
The purpose of this paper is to give some generalizations, in the context of Banach manifolds,  of  Sussmann's results about the orbits of families of vector fields (\cite{Su}). Essentially, we define the notion of "$l^1$-orbits" for any family of vector fields on a Banach manifold, and we prove, under appropriate assumptions,  that   such an orbit is a weak Banach submanifold.
  \end{abstract}

\bigskip



\section{Introduction}
 Let ${\cal X}$ be a family of local vector fields on a finite dimensional manifold $M$. According to the context of \cite{Su},  the orbit of ${\cal X}$ through $x\in M$ is the set of  points $\phi^{X_k}_{t_k}\circ\cdots\circ \phi^{X_1}_{t_1}(x)$ where $\{X_1,\cdots X_k\}$ is any finite  family of  vector fields in $\cal X$ and $\phi^{X_i}_{t}$ is the flow of $X_i$, $i=1,\cdots,k$. One most important  result of  H. Sussmann in $\cite{Su}$ is that each such an  orbit is an immersed submanifold of $M$. The proof of this result is founded on the two principal arguments
 \begin{enumerate}
 \item[(i)] enlargement of  ${\cal X}$ to the family $\hat{\cal X}$  of vector fields of type $(\phi^{X_p}_{t_p}\circ\cdots\circ\phi^{X_1}_{t_1})_*(X)$, for appropriate finite sets $\{X_1,\cdots, X_p,X\}\subset {\cal X}$ and  each orbit of $\hat{\cal X}$ is also an orbit of $\cal X$.
\item[(ii)] the distribution  $\hat{\cal D}$ generated by $\hat{\cal X}$ is integrable and each maximal integral manifold of $\hat{\cal D}$ is an orbit of $\hat{\cal X}$ and so also is an orbit of $\cal X$.
\end{enumerate}
As the dimension of $M$ is finite, the fundamental argument  for the proof of this last property   is that  $\hat{\cal D}$  is finite dimensional. \\

 For a generalization of such a result to Banach manifolds,  we can enlarge any family $\cal X$  in the same way as (i), but in (ii),  the argument of  finite dimension of the distribution $\hat{\cal D}$ is, of course, no more  valid. Naturally, we can hope that there exist some conditions under which analog  arguments  work for some "characteristic type" of families  of (local) vector fields on Banach manifolds.  So, given   a set ${\cal X}$ of local vector fields on a Banach manifold  $M$, after  having enlarged  ${\cal X}$  to a family  $\hat{\cal X}$  of vector fields  (in the same way  as (i)), we can look for  the orbits of $\hat{\cal X}$.  It is natural to consider the set of points of type
 \begin{eqnarray}\label{l1orbit}
 y=\phi^{X_n}_{t_n}\circ\cdots\circ \phi^{X_1}_{t_1}(x) \textrm{ or } y=\dis \lim_{k \ap \infty}\phi^{X_k}_{t_k}\circ\cdots\circ \phi^{X_1}_{t_1}(x)
\end{eqnarray}
 as an orbit through $x$ for any  finite or countable family $\{X_k\;, k \in A\}$ of vector fields in $\hat{\cal X}$. Note that,  if we restrict us  to finite sets $A$,   the binary relation  defined by\
 \begin{center}
 $y\backsim x$  if and only if    $y=\phi^{X_n}_{t_n}\circ\cdots\circ \phi^{X_1}_{t_1}(x)$,
 \end{center}
   is an equivalence relation. Moreover, in this case, there exists a piecewise smooth curve  which joins $x$ to $y$ and whose each connected part is tangent to  $X_i$ or $-X_i$ for some $i=1,\cdots, n$.

    Unfortunately, in the previous general case,  the associated  binary relation   clearly associated to (\ref{l1orbit}) is not any more a  relation of equivalence. The $\cal X$-orbit of $x$ will be the set of such points $y$ under some conditions so that  the associated binary relation is an equivalence relation.\\

 Given a family   $\xi\subset {\cal X}(M)$, a {\bf $\bf \xi$- piecewise smooth  curve} is a   piecewise smooth curve $\g:[a,b]\ap M$  such that each smooth part is tangent to $X$ or $-X$ for some $X\in \xi$.  In the context of (\ref{l1orbit}), for  such a point  $y$, there exists a family $\g_k:[0,T_k] \ap M$ of $\cal X$-piecewise smooth curves such that  the sequences of ends $x_k=\g_k(T_k)$ converge to $y$. When the sequence $T_k$ converges to some $T\in \R$ we have a continuous curve $\g:[0,T]\ap M$ such that $\g(0)=x$ and $\g(T)=y$. For  such a curve $\g$, there exists a countable  partition $t=(t_\a)_{\a\in A}$ of $[0,T]$ such that,  the restriction of $\g$ to $]t_\a,t_{\a+1}[$ is an integral curve of $X$ or $-X$, for some $X\in {\cal X}$. In particular the family $(\t_\a=t_{\a+1}-t_\a)_{\a\in A}$ belongs to $l^1(\N)$.  Such a curve will be called a {\bf $l^1$-curve of ${\cal X}$}.
The precise definition of an  orbit of ${\cal X}$   (see subsection \ref{pb}) is based on this notion of   $l^1$-curve but for the family $\hat{\cal X}$. Of course, we need some sufficient conditions under which $l^1$-curves exist. It is easy to see that condition of "local boundedness"  is a natural necessary condition, but, for the local existence, we need more: the local boundedness of the $s$-jets of vector fields of ${\cal X}$, for  sufficiently large $s>0$ (see subsection \ref{LB vect}). Under such assumptions, we can prove the existence of $l^1$-curves which are the integral curves of a vector field of type (see Theorem \ref{I}):
 $$Z(x,t,u)=\dis\sum_{\a \in A}u_\a(t)X_\a(x)$$
 where:

 $A$ is a finite,  countable or eventually uncountable set of indexes;

 $\xi=\{X_\a\}_{\a\in A}$ are defined on a same open set and  their $s$-jets are locally uniformly bounded (see Definition \ref{LBG});

$u=(u_\a)_{\a\in A} $ is a bounded integrable map from some interval $I$ to $l^1(A)$.\\

\noindent  In fact, in this context,  we get a {\bf flow $\Phi^\xi_u(t,\;)$} of such a vector field $Z$.\\

Let  $\xi=\{X_\a\}_{\a\in A}$ be a set  which satisfies a local boundedness condition for the $s$-jets for  sufficiently large $s>0$. The existence of $l^1$-curves which are integral curves of some $X\in \xi$ (or $-X$) on any subinterval $]t_\a,t_{\a+1}[$ associated to a countable partition of an interval $I$ is obtained by application of the previous result to $u=\G^\t=(\G^\t_\a)$ where $\G^\t_\a$ is the indicatrix function of  $]t_{\a},t_{\a+1}[$. Denote by $\Phi^\xi_\t(t,\;)$ the associated flow,  given any $x\in M$, for $T=||(t_{\a+1}-t_\a)_{\a\in A}||_1$,  $\t\ap \psi^x(\t)=\Phi^\xi_\t(T,x)$ is a map  from a neighborhood of $0\in l^1(A)$ into $M$ such that $\psi^x(0)=x$ and,  of class $C^{s-2}$, if the condition of local boundness of $s$-jets of elements of ${\cal X}$, are  satisfied (see Theorem 2).

Recall that our  purpose  is to prove, under appropriate assumptions, that each $\cal X$-orbit is a (weak) submanifold of $M$ as integral manifold of some distribution. According to the
proof of Sussmann's result, we first enlarge $\cal X$  into the set  $\hat{\cal X}$ given  by
$$\hat{\cal X}=\{Z=\Phi_*(\n Y),\; Y\in {\cal X},\; \Phi=\phi^{X_p}_{t_p}\circ\cdots\circ\phi^{X_1}_{t_1} \textrm{ for } X_1,\cdots, X_p\in {\cal X}  \textrm{  and appropriate }\n\in \R\} $$
(see subsection \ref{prelim}). From this set $\hat{\cal X}$, we associate  an appropriate pseudo-group ${\cal G}_{\cal X}$ of local diffeomorphisms, which is generated by flows of type $\phi^X_t$ with $X\in{\cal X}$ and  diffeomorphisms of type $\Phi_u^\xi(||\t||_1,.)$ (as we has seen previously) or its inverse for $\xi \subset\hat{\cal X}$. From this pseudogroup we get a  coherent and precise definition of  an {\bf orbit of $\bf{\cal X}$} or  {\bf ${\cal X}$ orbit} in short. Note that, under this definition, $\bf{\cal X}$ and $\bf{\hat{\cal X}}$ {\bf have the same orbits}, and moreover, if $y$ is in the orbit of $x$, there is a $l^1$ curve which joins $x$ to $y$ and whose smooth parts  are tangent to  vector fields of $\hat{\cal X}$. Note that the  binary relation associated to ${\cal G}_{\cal X}$  is then  an equivalence relation.  So, if $y$ belongs to the $\cal X$-orbit of $x$, either we have a $\cal X$- piecewise smooth curve which joins $x$ to $y$ or there exists a sequence $\g_k$ of $\cal X$-smooth piecewise curves whose origin is $x$ (for all curves) and whose sequence of ends converges to $y$
(see Proposition \ref{l1Xorbit} for a  complete description of a $\cal X$ orbit).\\

On the other hand, for any $x\in M$, under appropriate assumptions, we can associate the vector space $\hat{\cal D}_x=l^1(\hat{\cal X})_x$ which is  the set of all absolutely summable families $\dis\sum_{Y\in\hat{\cal X}} \t_Y Y(x)$.  In fact, in the same way, we can also associate the vector space  ${\cal D}_x=l^1({\cal X})_x$ generated by $\cal X$. Of course  ${\cal D}_x\subset \hat{\cal D}_x$ (see subsection \ref{characdist}) and we endox these vector spaces with a natural structure of Banach space. So we get weak distributions $\cal D$ and $\hat{\cal D}$ on $M$ such that $\hat{\cal D}$ is invariant  by any flow of vector fields in ${\cal X}$ and which is "minimal" for such a property (see Remark \ref{minSu}). Now, we need some conditions on $\hat{\cal X}$ which makes $\hat{\cal D}$ integrable. We will give two types of sufficient conditions.  \\

For the first one (called (H) in subsection \ref{loccond}), we assume that, for any $x\in M$, the Banach structure on $\hat{\cal D}_x$ is isomorphic to some $l^1(A)$ and  there exists a family $\{X_\a\}_{\a\in A}$ of vector fields defined around $x$, which are "locally uniformly bounded at order $s$" and such that $\{X_\a(x)\}_{\a\in A}$ is an unconditional symmetric basis of $\hat{\cal D}_x$. Under this assumption, $\hat{\cal D}$ is lower trivial (see subsection \ref{prelim}) but we cannot prove directly that $\hat{\cal D}$ is ${\cal X}^-_{\hat{\cal D}}$-invariant; in particular, {\bf we cannot use directly Theorem 1 of \cite{Pe}}. So we first prove that the map $\psi^x$, previously defined, gives rise to a local integral manifold of $\hat{\cal D}$ through $x$ of class $C^s$ for $s\geq 2$. This leads us to prove that  $\hat{\cal D}$ is ${\cal X}^-_{\hat{\cal D}}$-invariant and so we can now apply Theorem 1 of \cite{Pe}  and we finally get a smooth integral manifold of $\hat{\cal D}$. When $\hat{\cal D}$ is closed, we then obtain that each $l^1$-orbit has a structure of weak Banach manifold. {\it  Note that the assumption (H) is always satisfied when $\hat{\cal D}$ is finite dimensional (see Remark \ref{dimfin}). So this result can be seen as a generalization  of the proof Sussmann used in \cite{Su}}. \\

 The second sufficient conditions (called (H') in subsection \ref{stronginvol})  impose that $\hat{\cal D}$ is "upper trivial" (see section \ref{stronginvol}) and also some local involutivity conditions on  $\hat{\cal X}$. Under these conditions, by using a result of integrability from \cite{Pe}, we can show that $\hat{\cal D}$ is integrable and, when $\hat{\cal D}$ is closed,  each maximal integral manifold is a ${\cal X}$-orbit  (Theorem \ref{III"}). Moreover, if we consider   the family ${\cal X}^k$ defined by induction by:

${\cal X}^1={\cal X}$ and ${\cal X}^k={\cal X}^{k-1} \cup \{[X,Y],\; X\in {\cal X},\; Y\in {\cal X}^{k-1}\}$ for $k\geq 2$

we can associate, as previously, a weak distribution ${\cal D}^k=l^1({\cal X}^k)$. When such a distribution  satisfies the conditions (H') and is closed, we have ${\cal D}^k=\hat{\cal D}$ and so we get another sufficient conditions under which each ${\cal X}$-orbit is a weak manifold modelled on some $l^1(A)$. For the case where $\cal X$ is a finite family of global vector fields we get a new proof of the result of accessibility in \cite{Ro} (see Example \ref{serpent}). Moreover, when  $\cal X$ is a countable family of global vector fields, the reader can find an application of these results in \cite{PS}.\\

All these results can be naturally  applied in the context of control theory on Banach manifolds (Theorem \ref{V} and Theorem \ref{VI}).  These last   Theorems can be considered as a generalization of Sussmann's  accessibility results of \cite {Su} in finite dimension.\\

The paper is organized as follows.  In the next section, we study the problem of existence of $l^1$-curves. For any set $\cal X$ of vector fields which has the  "local boundedness of the $s$-jets of vector fields", we give sufficient conditions for the  existence of $l^1$-curves (Theorem 1) and we  apply this result to get $l^1$-curves tangent to $X\in{\cal X}$ or $-X$, on each subinterval associated to a countable partition. We also construct the map $\psi^x$ mentioned previously (Theorem 2).\\
 The notion of orbit of $\cal X$ or  $\cal X$-orbit, in short,  is  precisely defined in section 3. In the subsection \ref{prelim}, we construct  the announced enlargement $\hat{\cal X}$ of $\cal X$,  the associated  pseudogroup ${\cal G}_{\cal X}$ and we give a precise   definition of a $\cal X$-orbit. The following subsection is devoted to   all definitions and properties of distributions which will be used later.

Then the characteristic distributions  $\cal D$ and $\hat{\cal D}$ generated by ${\cal X}$ and $\hat{\cal X}$ respectively,  are defined in subsection \ref{characdist}. Finally, the main results of structure of weak Banach manifolds on ${\cal X}$-orbits are given and proved in section 4. In subsection \ref{loccond} under conditions (H) the corresponding result is given in Theorem 3. Under conditions (H'), the main results
 are given in Theorem 5.   Section 5   is devoted to some applications: on  one hand we obtain a new criterion of integrability of $l^1$-distributions  in Theorem \ref{IV} (see Remark \ref{intl1dis}). On the other hand,    we give  general results on  accessibility sets   as applications of the previous results on ${\cal X}$-orbits (Theorem \ref{V} and Theorem \ref{VI}).
 The last section is devoted to the proof of Theorem 2.\\

{\bf  Acknowlegements.} We would like to thank Professor  Tilmann Wurzbacher  for long and helpful discussions about  integrability of distributions and  for many useful  remarks about some proofs  in this work  and  also Professor Patrick Cabau for all his remarks and advices about  this paper.

\section{On $l^1$-integral curve of a uniformly locally bounded set of vector fields}
\subsection{Problem of existence of $l^1$-integral curve}\label{pb}
Let $M$ be  a smooth connected Banach manifold  modelled on a Banach space $E$.
 A {\bf local vector field } $X$ on  $M$ is a smooth section of the tangent bundle  $TM$   defined on an open set of  $M$ (denoted by  Dom($X$)). Denote by  ${\cal X}(M)$ the set of all local vector fields  on $M$. Such a vector field $X\in {\cal X}(M)$  has a flow  $\phi^X_t$ which is defined on a maximal  open set  $\O_X$ of  $\R\times M$. \\
 In this whole work, $A$,  $B$ and $\L$ will denote a finite  or a countable, eventually uncountable,  ordered  set of indexes. For such a countable set we shall often  identify this one with $\N$ as ordered set of indexes.

Consider a subset $\cal X$ of  ${\cal X}(M)$. As we have seen in the introduction, a curve $\g:[a,b] \ap M$ is called a $l^1$-integral curve of $\cal X$, if there exists a sequence  $t=(t_\a)_{\a\in A}$, where $A$ is a finite or countable set  such that:\\
-- $t_0=a$ and $t_{\a-1}\leq t_{\a}\leq b$ for $\a\in A$;\\
-- $t_n=b$ if $A$ is finite ($A\equiv \{1,\cdots n\}$) or $\dis\lim_{\a\ap \infty}t_\a=b$ (when $A$ is countable)\\
-- the restriction of $\g$ to each subinterval $]t_{\a-1},t_{\a}[$ is an integral curve of $X_\a$ or $-X_\a$ for some $X_\a\in {\cal X}$.\\
For such  a curve $\g$, the point $x_0=\g(a)$ (resp. $x_1=\g(b)$) is called  the origin (resp. the end) of $\g$ and we say that $x_0$ is joined to $x_1$ by a $l^1$-integral curve of $\cal X$.\\

It is clear that for any finite set $A=\{0,\cdots,n\}$ any $l^1$-integral curve is smooth by parts and, if we set  $\t_0=a$ and $\t_\a=t_\a-t_{\a-1}$ for $\a=1\cdots,n$, then there exist vector fields $X_1,\cdots X_n$ in $\cal X$ such that  for $\a=1,\cdots n$, we have:
    \begin{eqnarray}\label{compoflot}
\g(s)=\phi^{X_{\a}}_{s-t_{\a-1}}\circ\phi^{X_{\a-1}}_{\a-1}\circ\cdots\circ\phi^{X_1}_{\t_1}(\g(a)) \textrm{ for }  s\in [t_{\a-1},t_{\a}[ \textrm{ and } \a=1,\cdots,n
\end{eqnarray}
 Given a countable set $A\equiv \N$, and a $l^1$-integral curve $\g$  of $\cal X$, there exists a sequence of vector fields $\{X_\a\; , \a\in A\}$ in $\cal X$ such that (\ref{compoflot}) is true for all $\a\in A$. In particular, we must have:
\begin{eqnarray}\label{infinicompo}
\dis\lim_{\a\ap\infty}\phi^{X_{\a}}_{t_\a}\circ\cdots\circ\phi^{X_1}_{t_1}(\g(a)) =b.
\end{eqnarray}
We can cover such a  curve by a finite number of charts $(V_i,\phi_i)$, $i=1,\cdots r$ so that any $\g(]t_{\a-1},t_{\a}[)$ is contained in one domain $V_i$. Note that there exists one of these domains  which contains all $\g(]t_{\a-1},t_{\a}[)$ for $\a\geq \a_0$ for some $\a_0\in \N^*$ and we can assume that $V_r$ has this property. Now, on each $V_i$, a norm $||\,||_{\phi_i}$ can be defined on each fiber $T_xM$, for $x\in V_i$ by $||u||_{\phi_i}=||T_x\phi_i(u)||$ where $||\;||$ is a norm on $E$. From (\ref{infinicompo}), for any $\a\in A$, if  $\g(]t_{\a-1},t_{\a}[)\subset V_i$,  we must have
$$\sup \{||X_\a(\g(t))||_{\phi_i}\;, t\in ]t_{\a-1},t_{\a}[ \} \textrm{ is finite}$$

On the other hand, consider any countable  set  $A\equiv \N$  and any subset $\{X_\a\;, \a\in A\}$ of $\cal X$ such that Dom$(X_\a)$ contains $V$ and
$$\sup\{||X_\a(x)||_{\phi_i}\;, x\in V\;, \a\in A\} \textrm{ is finite }$$
Let be $\t=(\t_\a) \in l^1(A)$  such that $\t_\a>0$ for any $\a\in A$. Set $ t_0=0$  and $t_\a=\dis\sum_{i=1}^\a \t_i$ for $\a\in A$ and $T=\dis\lim_{\a\ap \infty}t_\a$. We set $\g(0)=x\in V$. If the flow $\phi^{X_1}_t(x)$ is defined for $t\geq \t_1$, we set $\g(t)=\phi^{X_1}(t)$ for $t\in [t_0,t_1]$. By induction,  suppose that we have defined $\g:[0,t_\a]\ap V$ such that $\g:[t_{i},t_{i+1}]\ap V$ is defined by $\g(t)=\phi^{X_i}_{t-t_i}(\g(t_{i}))$ for all $i=1\cdots \a$.  Then if  the flow $\phi^{X_{\a+1}}_t(\g(t_\a))$ is defined for $t\geq \t_{\a+1}$ then we put $\g(t)=\phi^{X_{\a+1}}_{t -t_\a}(\g(t_\a))$ for $t\in[t_\a,t_{\a+1}]$. So, when we can construct $\g$ at each step,  we get a $l^1$-integral curve of $\cal X$.
Consequently,  for the existence of $l^1$-integral curve associated to a countable  subset $\{X_\a\;, \a\in A\}$ of vector fields of $\cal X$,  we have to produce sufficient conditions under which  sequences of compositions  $$\phi^{X_{\a}}_{\t_\a}\circ\cdots \circ\phi^{X_{i}}_{\t_i}\circ\cdots\circ\phi^{X_1}_{\t_1},$$ converge when $\a\ap \infty$ and the limit defines a local diffeomorphism.  {\bf These conditions are assumptions  of uniform local boundness on the jets of vector fields} (see  next subsection).

\begin{rem}\label{homoborne}${}$\\
Consider a subset $\{X_\a\; \a\in A\}$ of $\cal X$ with the previous assumptions and $\t=(\t_\a) \in l^1(A)$. Recall that, for any local vector field $X$, for $\n\not= 0$, we have $\phi^X_t(x)=\phi^{X/\n}_{\n t}(x)$, when the second member is defined. It follows that given any $\n>0$
if a $l^1$-integral curve   $\g$  of $\{X_\a\; \a\in A\}$  is defined on $[0,T]$ as before, we can also define a $l^1$-integral curve $\bar{\g}$ of $\{\frac{1}{\n}X_\a, \a\in A\}$ in an obvious way on $[0,\n T]$ and we have $\bar{\g}(t))=\g(\n t)$ for any $t\in [0,T]$.
\end{rem}

\subsection{ Set of vector fields  uniformly locally bounded at order $s$}\label{LB vect}

Let $\Pi:TM\longrightarrow M$ be the tangent bundle of  $M$, with typical fiber $E$. Local vector fields on $M$ are local sections of this bundle.  Given  $X\in{\cal X}(M)$, the $s$-order  jet of $X$ at $x\in M$ is denoted by $J^s(X)(x)$. The set $J^s(TM)$ of  $s$-order  jets  of local vector fields on $M$ is a Banach bundle $\Pi^s:J^s(TM)\longrightarrow M$ of typical fiber  $E \times {\cal L}(E,E)\times
{\cal L}^2(E,E)\times\cdots,\times {\cal L}^s(E,E)$ where  ${\cal
L}^k(E,E)$, $2\leq k\leq s$  is the Banach space of symmetric $k$-linear maps from $E^k$ into  $E$ endowed with the usual norm (see for instance \cite{G} or  \cite{VE}). The typical fiber $E \times {\cal L}(E,E)\times {\cal L}^2(E,E)\times\cdots,\times
{\cal L}^s(E,E)$ of  $J^s(TM)$ is a Banach space for the norm  $||\,||_s$ which is the sum of the norm on $E$, the canonical norms on ${\cal L}(E,E)$  and on ${\cal
L}^k(E,E)$ for  $2\leq k\leq s$.\\

Consider  a chart  $(V, \phi)$ on  $M$ centered at  $x$.
On $V$ there exists
a trivialization $(\phi, \Phi)$ of  $[\Pi^s]^{-1}(V)$ on
$\phi(V)\times J^s(E)$ where  $J^s(E)=E \times {\cal L}(E,E)\times {\cal
L}^2(E,E)\times\cdots,\times {\cal L}^s(E,E)$ is the typical fiber.
On  $V$,  we have :
$$\forall y \in V, \; \Phi[J^s(X)(y)]=J^s(\phi_*X)(\phi(y))$$
So, on   $[\Pi^s]^{-1}(V)$, we have a norm  $|| \;||_\phi$
characterized by:
$$||J^s(X)(y)||_\phi= ||J^s(\phi_*X)(\phi(y))||_s$$

\begin{lem}\label{normeloc} (\cite{La})${}$\\
Let $V'$ be an open neighborhood of  $x$ having the same properties and $(\phi', \Phi')$the associated trivialization.
Denote by  $||J^s(X)(y)||_{\phi'}= ||J^s(\phi'_*X)(\phi'(y))||_s$
the associated norm on  $[\Pi^s]^{-1}(V')$. Then there  exists
a neighborhood $W\subset V\cap V'$  of $x$ and a constant  $C>0$
such that
$$\forall y \in W, \;||J^s[X](y)||_{\phi'}\leq C ||J^s[X](y)||_\phi$$
\end{lem}

\begin{defi}\label{LBG}${}$\\
Let  ${\cal X}$ be a set of local vector fields on $M$.
 Given $x\in M$, we say that
 $\cal X$ satisfies the condition
(LB(s))  at $x$ (Locally Bounded at order $s$), if there exist a chart
$(V_x, \phi)$ centered at $x$ and a constant $k>0$ such that:\\
for any $X\in {\cal X} $,  whose  domain
dom$(X)$ contains $V_x$, we have
\begin{eqnarray}\label{LBs}
\sup \{||J^s[X](y)||_\phi,\; X\in {\cal X} , \; y\in
V_x\}\leq k.
\end{eqnarray}
\end{defi}

\begin{rem}\label{remLB}${}$\\
 It follows from Lemma \ref{normeloc} that the property  (\ref{LBs}) does not depend neither on  the choice of the norm on $E$, nor on the choice  of the chart.
 \end{rem}

\begin{exes}\label{ex}${}$
\begin{enumerate}
\item[(1)] Let $E$ and $F$ be two Banach spaces and $T:F\ap E$ a continuous operator. Given any  finite or countable subset  $\{a_\a\;, \a\in A \}$   uniformly bounded  of  $F$ (i.e. $ || a_\a||\leq k $ for any $\a\in A$)   the assignment   $x \mapsto X_\a(x)=x+T(a_\a)$  is a vector field on $E$  and $\{X_\a\;,\a\in A \}$ satisfies the condition $LB(s)$ at any $x\in E$  and for any $s\in \N^*$.
\item[(2)] Let $L(F,E)$ be the set of continuous operators between the Banach spaces $F$ and $E$. Given a smooth map $\Phi: E\ap L(F,E)$, we denote by $\Phi_x$ the continuous operator associated to $x\in E$. By smoothness of $\Phi$, for any $x\in E$ and $s\in \N^*$, we can find an open neighborhood $U$ of $x\in E$ such that the jet of order $s$ of  $\Phi$ is bounded on $U$ (in the sense of Lemma \ref{normeloc}). Then, for any finite or countable subset  $\{a_\a\;,\a\in A\}$  uniformly bounded  of  $F$, denote by $X_\a$ the vector field on $E$ defined by $X_\a(x)=\Phi_x(a_\a)$. The set $\{X_\a\;,\a\in A \}$ satisfies the condition ($LBs$) at any $x\in E$  and for any $s\in \N^*$.
\item[(3)] Let ${\cal X}=\{X_1,\cdots,X_n\}$ be  a finite family of (global) vector fields on a Banach manifold $M$. Then ${\cal X}$ satisfies the condition (LBs), for any $s\in \N$.

\end{enumerate}
\end{exes}

\subsection{Sufficient conditions for the existence of $l^1$-integral curves}
\begin{nota}${}$\\
$\bullet$ $B(x,r)$ (resp. $B_f(x,r)$) denotes the open (resp. closed) ball centered at $x\in E$ of radius $r$ in the Banach space $E$.

\noindent $\bullet$  Given any Banach space $L$, if  $f:\R\times E\times L\ap E$ is a smooth map, we denote by $D_2f$ (resp. $D_3f$) the partial derivative relative to $E$ (resp. $L$).

\noindent $\bullet$ Let $\R^A$ will be the set of families $(u_\a)_{\a\in A}$  of  absolutely summable  real numbers where   $A$ is  countable or eventually uncountable set  of indexes  or   the set of finite real sequences $u=(u_1,\cdots u_n)$ if $A=\{1,\cdots, n\}$. We endow $\R^A$ with the norm
$$||u||_1=\dis\sum_{\a\in A}|u_\a|$$
It is well known that $(\R^A  ,||\;||_1)$ is a Banach space.

$\bullet$ Given any interval $J$ in $\R$ we denote by  $ L^1_b(J)$ the set of functions $u:J\ap R^A$ of class $L^1$ which are bounded. On $L^1_b(J) $ we define

\noindent - $\;\;\;||u||_{1} = \dis \int_{J} \dis\sum_{\a \in
A} |u_\a(t)|dt = \int_{J} \|u(t)\|_1dt$\hfill\\
- $\;\;\; ||u||_{\infty} = sup \{ \dis\sum_{\a \in A} |u_\a(t)|, \;
t \in J \} = sup\{||u(t) ||_1, \; t \in J \}$
\end{nota}

\noindent Given a finite, countable or uncountable ordered  set of indexes $A$, let $$\xi =
\{ X_\a, \; \a \in A, \; X_\a \in {\cal X}(M) \}$$
 be a set of vector fields on $M$ such that $\dis\bigcap_{\a\in A}$Dom$(X_\a)$ contains an open set $V$  of a chart $(V,
\phi)$ centered at $x$   such that the condition (LB(s+2)) at $x$ is satisfies
 for some  $s\in\N$. After restricting $V$ if necessary, we can suppose that  there exist $k>0$ such that
$$\sup \{||J^{s+2}(X)(y)||_\phi,\; X\in \xi , \; y\in
V\} \leq k$$

Without loss of generality, we can suppose that $V$ is an open set of the Banach space $E$.
To the previous set of vector fields  $\xi$, we can associate maps $Z$ of type:
\begin{displaymath}
\begin{tabular}{r c c l }
$Z:$&$J \times V \times L^1_b(J)$&$\longrightarrow$&$E$\\
&($t,x,u)$&$\longmapsto$&$Z(t,x,u)=\displaystyle\sum_{i \in J}
u_i(t)X_i(x)$
\end{tabular}
\end{displaymath}
It is easy to see that this map  $Z$  is of class
 $C^{s+1}$ relatively to the second variable.\\

 Given  such a map $Z$, let $J'$ be a subinterval of  $J$  and
$(t_0,x,u)\in J' \times V \times
L^1_b(I)$. A map  $f:J'\ap V$ is an {\bf integral curve} of $Z$, with initial condition $f(t_0)=x$ if
\begin{eqnarray}
\label{ci} \forall t \in I', \quad f(t)& = &x +
\int_{t_0}^{t}Z(s,f(s),u)ds
\end{eqnarray}

The following Theorem  gives the existence of a local flow  for $Z$:\\

\begin{theor}\label{I}${}$\\
Consider a fixed  $u$  in   $L^1_b(J)$, and we set   $c
=\|u\|_{\infty}$. Let  $(t_0, x_0, r,T',T_0)$ be an element of
$ J \times V \times {{\R}^*_+}^3 $ such that
$$]t_0 - T', t_0 + T'[ \subset J \; \; and  \;\;   B_f(x_0,2r)
\subset V$$
Moreover denote by: $$I_0=[t_0-T_0, t_0+T_0] \; \;
and  \; \;B_0 =B(x_0, r-kcT_0)$$
 If  $T_0 < min(\frac{r}{kc},T')$, then there exists a
flow $\Phi_u:I_0 \times B_0\ap V$, with the following properties:
\begin{enumerate}
\item for all $x$ in $B_0$,  each curve
$\Phi_u(.,x): I_0 \longrightarrow V$ is the unique   integral curve of
$Z$, with initial conditions $\Phi_u(t_0,x)=x$.
\item for all
$t\in I_0$, there exists an open connected  neighborhood  $U_0$ of $x_0$, contained in $B_0$ such that the map
$\Phi_u(t,.): U_0 \longrightarrow \Phi_u(t,U_0)$ is a
$C^s$-diffeomorphism.  Moreover, if  $D_2\Phi_u(t,.)$
and $D_2^2\Phi_u(t,.)$,denote the first and second derivative relative to the second variable,
we have :
\begin{displaymath}
\begin{tabular}{ l r c l }
$\forall x \in U_0,$&$D_2\Phi_u(t,x)$ & $=$& $ Id_E +\displaystyle \int_0^t D_2Z(s,\Phi_u(s,x),u)\circ D_2\Phi_u(s,x)ds$\\
&&&\\
&$D_2^2\Phi_u(t,x)$&$=$&$\displaystyle \int_0^t
(D_2^2Z(s,\Phi_u(s,x),u)\circ (D_2\Phi_u(s,x),D_2\Phi_u(s,x))
+$\\
&&&\\
&&&$\hfill \quad D_2Z(s,\Phi_u(s,x),u)\circ D_2^2\Phi_u(s,x))ds$\\
\end{tabular}
\end{displaymath}
\end{enumerate}
\end{theor}
\bigskip

This result is certainly well known for specialists. The reader can find a complete proof in \cite{La}.

Let  $\Phi$ and $\Psi$ be two local diffeomorphisms on $M$ which are defined on the domains  $\Omega_\Phi$ and
$\Omega_\Psi$ respectively. When $\Psi(\O_\Psi)\cap \O_\Phi \not=\emptyset$,
we can define the composition  $\Phi\circ\Psi$ which is a local diffeomorphism defined on $\Psi^{-1}[\Psi(\Omega_\Psi)\cap\O_\Phi]$. In this situation we will say that  $\Phi\circ\Psi$ is well defined. More generally,   we can
consider any finite composition $\Phi_n\circ\cdots\circ\Phi_1$ of local diffeomorphisms
$\Phi_1,\cdots,\Phi_n$ when successive compositions  $\Phi_i\circ(\Phi_{{i-1}}\cdots\circ\Phi_1)$ are well defined for $i=2,\cdots,n$. So, for a finite set  $A=\{1,\ldots n\}$,
and a finite set  $\xi=\{X_\a\}_{ \a\in A}$ of vector fields with the associate flows $\{\phi^{X_\a}_{t_\a}\}_{ \a\in A}$, it is clear that,  for $\tau=(\t_{1},\cdots,\t_{n})$,
under appropriate assumptions,  the composition
$\phi^\xi_\tau=\phi^{X_n}_{\t_{n}}\circ\ldots\circ\phi^{X_1}_{\t_{1}}$
is defined.  When $A$ is a countable or eventually uncountable ordered set of indexes we have the following result:

\begin{theor}\label{II}${}$\\
Let $\xi=\{X_\a\}_{ \a\in A}$ be a set of local vector fields such that Dom$(X_\a)$ contains $V$ for all $\a\in A$. Let be  $x_0\in V$ and  $r>0$ such that $ B_f(x_0,2r)$ is contained in $V$ and we assume that $\xi$ satisfies the condition (LB(s+2))  at $x_0$ where the relation  (\ref{LBs}) is true for all $y\in V$ and for the integer $s+2$.\\
Then, there exists an open connected neighborhood $U_0$ of  $x_0$, such that:
\begin{enumerate}
\item Fix any $\tau = (\t_\a)_{\a \in A}\in {\R}^A$ with $||\tau||_1\leq \displaystyle \frac{r}{k}$. Let $B$ be any countable subset of $A$ which contains all the indexes $\a$ such that $\t_\a\not =0$. Identifying the set   $B$ with $\N$ (as ordered sets), we denote by $\{\t_m,\; m\in B\}$ the associated subsequence of $\{\t_\a,\;\a\in A\}$. Then for any $x\in U_0$ we have:
\begin{enumerate} \item
$\phi^\xi_\tau(x) = \displaystyle \lim_{m \to
\infty}\phi^{X_m}_{\t_m} \circ \ldots \circ \phi^{X_1}_{\t_1}(x)$
exists. \item $\hat{\phi}^\xi_\tau(x) = \displaystyle \lim_{m \to
\infty} \phi^{X_1}_{-\t_1} \circ \ldots \circ \phi^{X_m}_{-\t_m}(x)$
exists. \item The map $\phi^\xi_\tau:x \longmapsto
\phi^\xi_\tau(x)$ is a local  $C^s$-diffeomorphism whose
inverse mapping  is $\hat{\phi}^\xi_\tau :x \longmapsto
\hat{\phi}^\xi_\tau(x)$.
\end{enumerate}
\item The map $\Psi^x$ defined in the following way:
\begin{displaymath}
\begin{tabular}{ r c c l r }
$\Psi^x : $&$\displaystyle B(0,\frac{r}{k}) $&$\longrightarrow$&$V$&\\
&$\tau$&$\longmapsto$&$\Psi^x(\tau)=\phi^\xi_\tau(x)$& is of class $C^s$.\\
\end{tabular}
\end{displaymath}
When the point $x$ will be fixed
we simply denote  $\Psi$ instead of $\Psi^x$.
\end{enumerate}
\end{theor}

The proof of this theorem  is long and technical, so it will be given in section \ref{appendice}.

\begin{rem}\label{inepB}${}$
\begin{enumerate}
\item  Denote by $\Phi^\xi_\t$ (resp. $\hat{\Phi}^\xi_\t$) the flow given in Theorem \ref{I} associated to $\xi$ and $u=\G^\t$ (resp. $\hat{u}=\hat{\G}^\t$) (see section \ref{appendice}). On the associated neighborhood $U$, we have

 $\hat{\Phi}^\xi_\t(t,z)=\Phi^\xi_\t(||\t||_1-t,\Phi^\xi_\t(-||\t||_1,z))$

 $\phi^\xi_\t(z)=\Phi^\xi_\t(||\t||_1,z)$

 $\hat{\phi}^\xi_\t(z)=[\phi^\xi_\t]^{-1}(z)=\hat{\Phi}^\xi_\t(||\t||_1,z)=\Phi^\xi_\t(-||\t||_1,z)$
 \item In fact,  both limits $\phi^\xi_\t(x)$ and $\hat{ \phi}^\xi_\t(x)$  do not depend on the choice of the set $B$ but only depend on the countable set  $A_\t=\{\a\in A \textrm{ such that } \t_\a\not=0\}$. Moreover, the set $A_\t$ is independent of $x\in U_0$.
\item To each $\t$  the  associated  set $A_\t=\{\a\in A \textrm{ such that } \t_\a\not=0\}$ can be written $A_\t= \{\a_k,k\in \N\}$ or $A_\t =\{\a_k\;, k=1\cdots, n\}$. Consider the associated    subdivision $\{t_{\a_k}\}_{k\in \N}$ of the interval $[0,T]$ defined by:\\
$t_0=0\leq t_1=|\t_{\a_1}|\leq\cdots\leq t_i=\dis\sum_{k=1}^i|\t_{\a_k}|\leq \cdots \leq T=\dis\sum_{\a\in A_\t}|\t_{\a}|$.\\
 Fix some $x\in U_0$ and  let $(x_k)_{\a_k \in A_\t}$ be the sequence defined by:

$x_0=x$, and for $\a_k\in A_\t$,  $x_k=\phi^{X_{\a_k}}_{\t_{\a_k}}(x_{k-1})=\phi^\xi_\t(x_{k-1})$. \\Then the curve $\g:[0,T]\ap M$ defined by $\g(s)=\phi^{X_{\a_k}}_{s-t_{k-1}}(x_{k-1})=\Phi^\xi_\t(s,x)$ for $s\in[t_{k-1},t_k[$ is a $l^1$-curve which joins $x$ to $\Psi^\xi_\t(x)$. On the other hand, to $\hat{ \phi}^\xi_\t$ we can associate the curve $\hat{\g}:[0,T]\ap M$  defined by $\hat{\g}(s)=\g(T-s)$. So $\hat{\g}$ joins $\g(0)=\phi^\xi_\t(x)$ to $x$. We also call such a curve,  the $l^1$-curve associated to $\hat{ \phi}^\xi_\t$.
\end{enumerate}
\end{rem}

\section{The orbits   of ${\cal X}$ or  $\cal X$-orbits}
\subsection{ Definition of  an orbit of ${\cal X}$}\label{prelim}

In this section, we consider a fixed set ${\cal X}$  of vector fields on $M$  with the following properties:
\begin{enumerate}
\item[(Hi)] $M= \dis\bigcup_{x\in {\cal X}}$ Dom$(X)$
\item[(Hii)] there exists $s \geq 0$ with the following property:
for any $x\in M$ there exists a chart  $(V_x,\phi)$ centered at $x$ such that for  the set ${\cal X}_{x}$ of vector fields $X\in {\cal X}$ whose  Dom$(X)$ contains $x$ we have
$$\sup \{||J^{s+2}[X](x)||_\phi,\; X\in {\cal X} \}<\infty.$$
\end{enumerate}

The announced  enlargement ${\hat{\cal X}}$ of $\cal X$ is  obtained  from the following Lemma:

 \begin{lem}\label{htX}${}$\\
Let $(V_x,\phi)$ be a chart  centered at $x$  and a constant $k$ such that
$$\sup \{||J^{s+2}[X](x)||_\phi,\; X\in {\cal X} \}\leq k.$$
Let $\hat{\cal X}_x$ be the set of local vector fields  of type $Y=(\phi^{X_p}_{t_p}\circ\cdots\circ\phi^{X_1}_{t_1})_*(\n.X)$, for any $\n>0$,  where $X_1,\cdots,X_p,X$ belongs to $\cal X$, whose domain contains $x$ and such that
\begin{eqnarray}\label{2LBs+2}
||J^{s+2}[Y](x)||_\phi \leq k
\end{eqnarray}
 We set
$$\hat{\cal X}=\dis\bigcup_{x\in M}\hat{\cal X}_x$$
\begin{enumerate}
\item[(i)] $\hat{\cal X}$ contains ${\cal X}$ and satisfies the conditions ${\textrm(Hi)}$ and ${\textrm(Hii)}$.
\item[(ii)]Let   $\hat{\hat{\cal X}}$ be  the set of vector fields obtained from $\hat{\cal X}$  in the same way as $\hat{\cal X}$ from ${\cal X}$. Then,
we have $\hat{\hat{\cal X}}=\hat{X}$
 \end{enumerate}
 \end{lem}
\bigskip

\begin{rem}\label{hatflow}${}$\\
According to Remark \ref{homoborne}, the flow of any vector field $Y=(\phi^{X_p}_{t_p}\circ\cdots\circ\phi^{X_1}_{t_1})_*(\n.X)$  can be written
\begin{eqnarray}\label{Yflow}
\phi^Y_\t=\phi^{X_1}_{-t_1}\circ\cdots\circ\phi^{X_p}_{-t_p}\circ\phi^X_{ \t/\n}\circ \phi^{X_p}_{t_p}\circ\cdots\circ\phi^{X_1}_{t_1}
\end{eqnarray}
\end{rem}
\bigskip
\noindent\begin{proof}\so{\it proof}${}$\\
\noindent Let be $x\in M$. For any $X\in {\cal X}_{x}$, we have $(\phi^X_0)_*X=X$ and, by construction, the  vector fields in $\hat{\cal X}_{x}$ satisfies the condition (\ref{2LBs+2})   with the same constant $k$, so ${\cal X}_{x}$ is contained in  $\hat{\cal X}_{x}$ . If follows that $\hat{\cal X}$ satisfies (Hi). The condition  (Hii) follows from the definition of $\hat{\cal X}_x$.\\

By construction,  ${\hat{\hat{\cal X}}}_x$ is the set of vector fields $Z=(\phi^{Y_p}_{t_p}\circ\cdots\circ\phi^{Y_1}_{t_1})_*(\n.Y)$  where $Y_1,\cdots,Y_p,Y$ belongs to $\hat{\cal X}$ for some $\n>0$. As  we have $Y=(\phi^{X_q}_{t_q}\circ\cdots\circ\phi^{X_1}_{t_1})_*(\n'.X)$,    from  Remark \ref{hatflow} we get
$$Z=(\phi^{X_m}_{s_m}\circ\cdots\circ\phi^{X_1}_{s_1})_*(\n \n' X)$$
 for appropriate vector fields $X_1,\cdots, X_m,X$ in $\cal X$ and appropriate real  values $s_1,\cdots, s_m$.

Now, on the considered chart $(V_x,\phi)$, we have $||J^{s+2}[Y](x)||_\phi  \leq k, X\in {\cal X}$. So we also have  $||J^{s+2}[Z](x)||  \leq k$. We conclude that $Z$ belongs to $\hat{\cal X}$

\end{proof}\\

Let  {\bf $\bf{\cal G}_ {\cal X}$ be the pseudogroup} of local diffeomorphisms $\Psi$  which are  defined in the following way: \\

$\Psi= \phi_n\circ  \cdots\circ\phi_k\circ\cdots\circ\phi_1$ when these compositions are well defined and \\
where $\phi_k$ is a local diffeomorphism of one of the following type

(i)$\phi^X_{\t_k}$ for some $X\in {\cal X}$ and $\t_k\in \R$

(ii) $\phi^{\xi_k}_{\tau_k}$ or $[\phi^{\xi_k}_{\tau_k}]^{-1}$  as defined in Theorem \ref{II},
where   $\xi_k=\{X_\a \;, \a\in A_k\}$  is a finite or countable

 subset of $\hat{\cal X}$
and $\tau_k\in \R^{A_k}$.\\

\begin{com}\label{Xl1curve}${}$
\begin{enumerate}
 \item From (\ref{Yflow}) any  flow $\phi^Y_\t$ for $Y\in \hat{\cal X}$ belongs to ${\cal G}_{\cal X}$.
 \item Let be $\Psi= \phi_n\circ\cdots\circ\phi_1\in {\cal G}_{\cal X}$. By construction of $\Psi$,  to each $\phi_k$ is associated  a family $\xi_k=\{X_\a \;, \a\in A_k\}$  which is a finite or countable subset of $\hat{\cal X}$
and $\tau_k\in \R^{A_k}$,  we have a real positive number  $\dis\sum_{k=1}^n||\t_k||_1<\infty$ associated to $\Psi$. If $\phi_k$ is of type (ii),  according to Remark \ref{inepB} {\it 1.}, denote by $\Phi_k$ the flow associated to each $\xi_k$ with $u=\G^{\t_k}$ or $u=\hat{\G}^{\t_k}$ if $\phi_k=\phi^{\xi_k}_{\t_k}$ or $\phi_k=[\phi^{\xi_k}_{\t_k}]^{-1}$ respectively. If $\phi_k$ is of type (i)  $\xi_k$ is reduced to some $X_k\in {\cal X}$ and we have $\Phi_k(t,y)=\phi^{X_k}_t(y)$.

Take  any pair $(x,y)\in M^2$ such that $y=\Psi(x)$.  We set $t_0=0$ and $t_k=\dis\sum_{i=1}^k||\t_i||_1$ for $k=1,\cdots,n$. Consider the sequence $(x_k)$ defined by $x_0=x$ and $x_k=\Phi_k(\t_k,x_{k-1})$. So for each $k$, we can consider the $l^1$-curve $\g_k:[t_{k-1},t_k]\ap M$ defined by $\g_k(t)=\Phi_k(t-t_{k-1},x_{k-1})$ (see Remark \ref{inepB} {\it 3.}). By construction, we have $\g(t_k)=x_k$ and  $y=x_n$. So if $T=\dis\sum_{k=1}^n ||\t_k||_1$ we get a  sequence of  $l^1$-curve $\g=[0,T] \ap M$,  defined by $\g_{|[t_{k-1},t_k[}=\g_k$,  such that $\g(0)=x$ and $\g(T)=y$.
\item Given a family   $\xi\subset {\cal X}(M)$, recall that  a {\bf $\bf \xi$- piecewise smooth curve} is a  piecewise smooth curve $\g:[a,b]\ap M$  such that each smooth part is tangent to $X$ or $-X$ for some $X\in \xi$. \\

When  $y=\phi^Y_\t(x)$ for $Y\in \hat{\cal X}$, from (\ref{Yflow}),  we can clearly  associate a $\xi$- piecewise smooth curve which joins $x$ to $y$.

 Now,  consider any  $\xi=\{X_\a,\; \a\in A\}\subset\hat{\cal X}$ and $\t$ small enough such that  $\phi^\xi_\t$ is defined and consider $y=\phi^\xi_\t(x)$. If $A=\{1,\cdots,n\}$ is finite, from the previous argument, there exists a family $\xi_n\subset {\cal X}$  and an associate $\xi_n$- piecewise smooth curve $\g'_n$ which joins $x$ to $y$. On the other hand, if $A$ is countable, to each $k\in A$, we can associate a family $\xi_k\subset{\cal X}$  and a $\xi_k$- piecewise smooth curve $\g'_k$ which joins $x=x_0$ to $x_k$ (as defined in Remark \ref{inepB} {\it 3.}). So we get a sequence of  ${\cal X}$- piecewise smooth curves whose   origin is $x_0$ (for all curves) and  whose sequence of ends converges to $y$. Note that, for Theorem \ref{II}, the same result is true for any pair $(z,\Phi^\xi_\t(z))$ for any $z$ in some neighbourhood $U$ of $x$ we have
 \begin{eqnarray}\label{phixi}
 {\phi}^\xi_\tau(z) = \displaystyle \lim_{m \to\infty} \phi^{X_m}_{-\t_m} \circ \ldots \circ \phi^{X_1}_{-\t_1}(z) \textrm{ for any } z\in U
\end{eqnarray}
where $\xi=\{X_k,\; k\in A\}\subset\hat{\cal X}$ and $\t=(t_k)_{k\in A}$.

 From (\ref{Yflow}) to each  finite sequence  $\phi^{X_m}_{-\t_m} \circ \ldots \circ \phi^{X_1}_{-\t_1}(z)$, we can associate a $\cal X$- piecewise smooth curve $\g'_m$ which joins $z$ to $z_m= \phi^{X_m}_{-\t_m} \circ \ldots \circ \phi^{X_1}_{-\t_1}(z)$. So  given $ {\phi}^\xi_\t$, for any $z\in V$ we have a family of $\cal X$- piecewise smooth curves $\g'_m$ whose origin is $z$ and whose sequence of ends converges to ${\phi}^\xi_\tau(z)$

 Now, consider the case $y=\hat{\phi}^\xi_\t(x)=[\phi^\xi_\t]^{-1}(x)$. Again, from Theorem  \ref{II}, there exists some open neighbourhood $U$ of $x$ such that $ \hat{\phi}^\xi_\t$ is a local diffeomorphism on  $U$ of $x$ and we have
\begin{eqnarray}\label{hatphixi}
\hat{\phi}^\xi_\tau(z) = \displaystyle \lim_{m \to\infty} \phi^{X_1}_{-\t_1} \circ \ldots \circ \phi^{X_m}_{-\t_m}(z) \textrm{ for any } z\in U
\end{eqnarray}
where $\xi=\{X_k,\; k\in A\}\subset\hat{\cal X}$ and $\t=(t_k)_{k\in A}$.

Again from (\ref{Yflow}) to each  finite sequence  $\phi^{X_1}_{-\t_1} \circ \ldots \circ \phi^{X_m}_{-\t_m}(z)$, we can associate a $\cal X$- piecewise smooth  curve $\g'_m$ which joins $z$ to $z_m= \phi^{X_1}_{-\t_1} \circ \ldots \circ \phi^{X_m}_{-\t_m}(z)$. So  given $ \hat{\phi}^\xi_\t$, for any $z\in V$ we have a family of $\cal X$- piecewise smooth curves $\g'_m$ whose origin is $z$ and whose sequence of ends converges to $\hat{\phi}^\xi_\tau(z)$. This is in particular true for the previous fixed pair $(x,y)$.

 ${}\;\;\;\;$ In the general case when $y=\Phi(x)$ for some $\Phi\in {\cal G}_{\cal X}$, we have $\Phi= \phi_n\circ\cdots\circ\phi_1\in {\cal G}_{\cal X}$. Set $x_1=\phi_1(x)$. From the previous partial results, either we have a ${\cal X}$- piecewise smooth curve which joins $x$ to $x_1$ or there exists a sequence $\g_k$ of $\cal X$- piecewise smooth curves whose origin is $x$ (for all curves) and whose sequence of ends converges to $x_1$.  At first, assume that we are in the first case. Now, if we have a ${\cal X}$- piecewise smooth curve which joins $x$ to $x_1$, applying the previous argument in $x_1$ by concatenation, we get either a ${\cal X}$- piecewise smooth curve which joint $x$ to $x_2=\phi_2(x_1) $ or we get a sequence of a sequence  of $\cal X$- piecewise smooth curves whose origin is $x$ (for all curves) and whose sequence of ends converges to $x_2$. If we are in the second case, where $V$ is a neighborhood  of  $x_1$ on which (\ref{phixi}) or (\ref{hatphixi}) is true.  For $k$ large enough, $\g_k(x_1)$ belongs to $V$. So for each $k$, we have a family of $\cal X$- piecewise smooth curves $\g'_{(k,n)}$  whose origin is $\g_k(x)$ (for all curves) and whose sequence of ends converges to $\phi_2(\g_k(x_k))$. As $\dis\lim\phi_2(\g_k(x_k))=x_2$, there exists an  increasing sequence $n_k$ such that  the sequence $\hat{\g}_k$ of the concatenations $\g_k$ with $\g'_{(k,n_k)}$ is a sequence of $\cal X$- piecewise smooth curves  whose origin is $x$ (for all curves) and whose sequence of the  ends converges to $x_2$. By finite induction on $k$, we get the same result for the pair $(x,y)$.
   \end{enumerate}
\end{com}

\noindent  To ${\cal G}_{\cal X}$ is naturally associated the following equivalence relation
on  $M$:
$$x\equiv y \textrm{ if and only if there exists } \Phi\in{\cal G}_{\cal X} \textrm{ such that } \Phi(x)=y$$

\noindent {\bf An equivalence class is called a $ l^1$-orbit of $\bf{\cal X}$ or a $\bf {\cal X}$-orbit}.\\

\noindent The term   "$ l^1$-orbit" will be justified by the following  result which sums up the previous commentaries  and Lemma \ref{htX} part (ii):
\begin{prop}\label{l1Xorbit}${}$
\begin{enumerate}
\item Each point of the $ {\cal X}$-orbit of $x$ can be joined from $x$ by a $l^1$-curve whose each connected  smooth part is tangent to  $Y$ or $-Y$ for some $Y\in {\hat{\cal X}}$.
\item  For each pair $(x,y)$ in the same $ {\cal X}$-orbit, either we have a ${\cal X}$-piecewise smooth curve which joins $x$ to $y$ or there exists a sequence $\g_k$ of $\cal X$- piecewise smooth curves whose origin is $x$ (for all curves) and whose sequence of the  ends converges to $y$.
\item Let ${\cal G}_{\hat{\cal X}}$ be the pseudogroup naturally associated  to ${\hat{\cal X}}$. Then we have ${\cal G}_{\hat{\cal X}}={\cal G}_{\cal X}$. In particular each $\hat {\cal X}$-orbit is a $ {\cal X}$-orbit.
\end{enumerate}
\end{prop}

\subsection{Preliminaries on weak distributions}\label{weakdistrib}

${}\;\;$ {\it Recall that,  according to the proof of Sussmann's  theorem on reachable sets in \cite{Su},   we want to associate to $\cal X$ and $\hat{\cal X}$  weak distributions ${\cal D}$ and $\hat{\cal D}$ respectively, such that ${\cal D}_x\subset \hat{\cal D}_x$ for any $x\in M$, $\hat{\cal D}$ is invariant by any flow of vector fields in $\cal X$
  and which is minimal (in some sense) for these properties.\\
  Before  beginning this construction, we need to recall some definitions on weak distributions  which will be used   in the next subsection.}\\

$\bullet$ Given a finite or countable or eventually uncountable ordered set  $A$ of indexes,  a family $\{\e_\a,\; \a\in A\}$ is said to be an unconditional basis of $\R^A$ if, for every $\t \in \R^A$  there is a unique family of scalars $\{\t_\a\,;\a\in A\}$ such that $\t=\dis\sum_{\a\in A}\t_\a \e_\a$ (unconditional convergence); such a basis is symmetric if for any sequence $(\a_k) \in A$ with $k\in K\subset \N$, the basic sequence $\{\t_{\a_k}\;, k\in K \}$ is equivalent to the canonical basis of $\R^K$ (see for instance \cite{LT}).  It is well known that all unconditional symmetric basis of $\R^A$ are  equivalent to the canonical basis of $\R^A$.\\

$\bullet$ A {\bf weak submanifold}  of $M$ of class $C^p$ (resp. smooth) is a pair $(N,f)$ of a connected  Banach manifold  $N$ of class $C^p$ (resp. smooth) (modeled on a Banach space $F$) and a map  $f:N\ap M$  of class $C^p$ (resp. smooth) such that: (\cite{El},\cite{Pe})

-- there exists a continuous injective linear  map $i:F\ap E$ between these two Banach spaces

-- $f $ is injective and the tangent map $T_xf:T_xN\ap T_{f(x)}M$ is injective for all $x\in N$.\\

 Note that for  a weak submanifold $f:N\ap M$, on the subset $f(N)$ in $M$ we have two topologies:

--  the induced topology from $M$;

--  the topology for which $f$ is a homeomorphism from $N$ to $f(N)$.\\
With this last topology, via $f$, we get a structure of Banach manifold modeled on $F$.
 Moreover, the inclusion from $f(N)$ into $M$ is continuous as a map from the Banach manifold $f(N)$ to $M$.
 In particular, if $U$ is an open set of $M$, then, $f(N)\cap U$ is an open set for the topology of the Banach manifold  on $f(N)$.\\

$\bullet$ According to \cite{Pe}, a  {\bf weak distribution} on a $M$  is an assignment  $ {\cal D}: x \mapsto {\cal D}_x$  which, to every $x\in M$, associates  a vector subspace  ${\cal D}_x$ in $T_xM$  (not necessarily closed) endowed with a norm $||\;||_x$ such that   $({\cal D}_x, ||\;||_x)$ is a Banach space (denoted by $\tilde{\cal D}_x$) and such  that the natural inclusion $i_x : \tilde{\cal D}_x \ap T_xM$ is continuous.\\

When ${\cal D}_x$ is closed,  we have a natural Banach structure  on $\tilde{\cal D}_x$, induced by the Banach structure on $T_xM$, and so we get the classical definition of a distribution;  in this case we will say that  $\cal D$ is {\bf closed}.

\noindent A vector field $Z\in{\cal X}(M)$ is {\bf tangent} to $\cal D$, if for  all $x\in$ Dom($Z$), $Z(x)$ belongs to ${\cal D}_x$. The set of local vector fields tangent to $\cal D$  will be denoted by {\bf $\bf {\cal X}_{\cal D}$}.

$\bullet$ We say that  $\cal D$ is {\bf is generated by a subset  }  ${\cal X}\subset {\cal X}(M)$ if, for every $x\in M$, the vector space ${\cal D}_x$ is the  linear hull  of the set $\{Y(x)\;,\;Y\in{\cal X}\;,\;x\in$ Dom$(Y)\}$.  \\

For a weak distribution $\cal D$, on $M$ we have the following definitions:\\

$\bullet$   $\cal D$ is   {\bf lower (locally) trivial} at $x$ if there exists an open neighborhood $V$ of $x$ in $M$,  a smooth map $\Phi:\tilde{\cal D}_x\times  V \ap TM$  (called {\bf  lower trivialization}) such that :
\begin{enumerate}
\item[(i)]  $ \Phi(\tilde{\cal D}_x\times\{y\})\subset {\cal D}_y$ for each $y\in V$

\item[(ii)] for each $y\in V$,  $\Phi_y\equiv \Phi(\;,y):\tilde{\cal D}_x\ap T_yM$ is a continuous operator  and $\Phi_x:\tilde{\cal D}_x\ap T_xM$  is the natural inclusion $i_x$

\item [(iii)] there exists a  continuous operator $\tilde{\Phi}_y: \tilde{\cal D}_x\ap \tilde{\cal D}_y$ such that $i_y\circ \tilde{\Phi}_y=\Phi_y$, $\tilde{\Phi}_y$ is an isomorphism from $\tilde{\cal D}_x$ onto ${\Phi}_y(\tilde{\cal D}_x)$
and  $\tilde{\Phi}_x$ is the identity of $\tilde{\cal D}_x$.
\end{enumerate}
We say that  $\cal D$ is  {\bf lower (locally) trivial}  if it is lower trivial at any $x\in M$.\\

$\bullet $  ${\cal D}$ is called a {\bf $ { \bf l^1}$- distribution}  if   each  Banach space  $\tilde{\cal D}_x$ is isomorphic to  $\R^A$,  for some appropriate finite, countable or eventually uncountable ordered  set $A$ of indexes (which depends of $x$).\\

$\bullet$ an {\bf integral manifold} of class $C^p$, with $p\geq 1$, (resp. smooth) of $\cal D$  through $x$ is a weak  submanifold $f:N\ap M$  of class $C^p$ (resp. smooth) such that there exists $x_0\in N$ with $f(x_0)=x$ and $T_zf(T_zN)={\cal D}_{f(z)}$ for all $z\in N$. An integral manifold  through $x\in M$ is called {\bf maximal} if, for any integral manifold $g:L\ap M$ through $x$, the set $g(L)$ is an open submanifold of $f(N)$, according to the  structure of Banach manifold on  $f(N)$ induced by $N$ via $f$.\\

$\bullet$   $\cal D$ is called {\bf integrable} of class $C^p$ (resp. smooth) if for any $x\in M$ there exists an integral manifold $N$ of class $C^p$ (resp. smooth) of $\cal D$ through $x$.\\

$\bullet$ if  $\cal D$ is generated by ${\cal X}\subset {\cal X}(M)$, then    $\cal D$ is  called {\bf ${\bf {\cal X}}$- invariant} if for any   $X\in {\cal X}$, the tangent map $T_x\phi^X_t $ send ${\cal D}_x$ onto $ {\cal D}_{\phi^X_t(x)}$ for all $(t,x)\in \O_X$. ${\cal D}$ is {\bf invariant} if ${\cal D}$ is $ {\cal X}_{\cal D}-$ invariant. \\

\subsection{Characteristic distribution associated to $\cal X$}\label{characdist}


\noindent Consider any set $\cal Y$ of local vector fields such that,   conditions (Hi) and (Hii) are satisfied. We denote by ${\cal Y}_x$ the set of vector fields $Y\in {\cal Y}$ such that $x$ belongs to Dom$(Y)$.  The distribution $l^1({\cal Y})$ defined by:
$$l^1({\cal Y})_x=\{X\in T_xM \textrm{ such that } X=\dis\sum_{Y\in{\cal Y}_x} \t_YY(x) \textrm { with } \dis\sum_{Y\in{\cal Y}} |\t_Y| \textrm{ summable }\}$$
is called {\bf the  $l^1$-characteristic distribution} generated by $\cal Y$.

 For $x\in M$ fixed, let  $\L$ be any (ordered) set of indexes of same cardinal as ${\cal Y}_x$ so that each element of ${\cal Y}_x$ can be indexed as $Y_\l,\; \l\in \L$. We then have a surjective linear map:
$T:l^1(\L)\ap T_xM$ defined by $T((\t_\l)_{\l\in \L})=\dis\sum_{\l\in \L} \t_\l Y_\l(x) $ and whose range is $ l^1({\cal Y})_x$. So we get a bijective continuous  map $\bar{T}$ from the quotient $l^1(A)/\ker T$ onto $l^1({\cal Y})_x$. So we can put on  $l^1({\cal Y})_x$  a structure of Banach space such that  $\bar{T}$  is an isometry. Finally, $l^1({\cal Y})$ is a weak distribution. $l^1({\cal Y})_x$ will always be equipped with this  Banach structure.\\

\begin{rem}\label{compl1}${}$
\begin{enumerate}
\item For the the existence of  $l^1({\cal Y})_x$ we only need  that  for all $X\in{\cal Y}_x$
$$\sup \{||X(x)||_\phi,\; X\in {\cal X} \}<\infty.$$
So the condition (Hii) is much too strong  in this way. However,  independently of the existence of  $l^1({\cal Y})_x$, in this paper, we need to consider the set  ${\cal Y}$ of local vector fields which satisfies condition (Hii).
\item The Banach space $l^1({\cal Y})_x$ is isomorphic to $l^1(A)$ for some ordered set $A$ if and only if, with the previous notations, $\ker T$ is complemented. In this case, $A$ has the same cardinal as $\L$ (see \cite{Ko}). In particular, if the distribution $l^1({\cal Y})$ is  upper trivial (see subsection \ref{stronginvol}), then   $l^1({\cal Y})_x$ is isomorphic to some $l^1(A)$ for any $x\in M.$
\end{enumerate}
\end{rem}
The {\bf characteristic  distribution  ${\bf \cal D}$ associated to} $\cal X$ is defined by:
$${\cal D}_x=l^1({\cal X}_x)$$
 Note that,  from assumptions (Hi) and  (Hii), ${\cal D}_x$ is well defined for any $x\in M$. Moreover, the natural inclusion of ${\cal D}_x$ into $T_xM$ is continuous.\\
 In the same way,   the {\bf characteristic distribution ${\bf \hat{\cal D}}$  associated to} $\hat{\cal X}$ is defined by:

 $$\hat{\cal D}_x= l^1(\hat{\cal X }_x)$$
From Lemma \ref{htX} part(i) it follows that $\hat{\cal D}$ is well defined  and, again,  the natural inclusion of $\hat{\cal D}_x$ in $T_xM$ is continuous. Moreover, as ${\cal X}_x\subset \hat{\cal X}_x$, we have ${\cal D}_x\subset \hat{\cal D}_x$ for any $x\in M$. The other relative properties of $\cal D$ and $\hat{\cal D}$ are given in the following  Proposition.

\bigskip
\begin{prop}\label{ProI} ${}$
\begin{enumerate}
\item
 $\hat{\cal D}$ is ${\cal X}$-invariant and also $\hat{\cal X}$-invariant.
 \item Let ${\cal Y}$ be any family of local vector fields which satisfies (Hi) and (Hii) and which contains ${\cal X}$. If the associated distribution $l^1({\cal Y})$ is ${\cal X}$-invariant then $l^1({\cal Y})_x$ contains $\hat{\cal D}_x$ for any $x\in M$. In particular, if $\cal D$ is ${\cal X}$-invariant, then ${\cal D}=\hat{\cal D}$.
 \item Given  $x\in M$ and assume that we have the following properties:
\begin{enumerate}
\item[(i)]  there exists  a finite countable or  eventually uncountable set $A$ of indexes such that $\hat{\cal D}_x$ is isomorphic to $\R^A$
\item[(ii)] there exists  a chart domain $V_x$ centered at $x$ and a family $\{X_\a, \a\in A\}\subset \hat{\cal X}_{x}$
such that $\{X_\a, \a\in A\}$ satisfies the condition (LB(s+2)) on $V_x$, for some $s>0$,  and,
$\{X_\a(x), \a\in A\}$ is a symmetric unconditional basis of $\hat{\cal D}_x\equiv \R^A$.
\end{enumerate}
 Then, there  exists a weak Banach manifold $\Theta:
B(0,\r)\ap M$  of class $C^s$, which is an integral manifold of $\hat {\cal D}$ through $x$,
where $B(0,\r)$ is the open ball in the Banach
space $\R^A$. Such a manifold will be called a {\bf slice centered at $x$}.
 \item Let $f:N\ap M$ be a smooth connected integral manifold such that $x\in f(N)$.  Assume that the hypothesis of part {\it 3} are satisfied at $x$. Then, for $\r$ small enough
, $\Theta(B(0,\r))$ is contained in  $f(N)$ and $f^{-1}(\Theta(B(0,\r)))$ is an open set in $N$.
\end{enumerate}
\end{prop}

\begin{rem}\label{minSu}${}$\\
Classically,  a distribution on $M$ is an assignment $\D: x\mapsto \D_x$ where $\D_x$ is a vector subspace of $T_xM$. As in \cite{Su}, on the set of distributions, we can consider the partial order:

$\D\subset \D'$ if and only if $\D_x\subset \D'_x$ for any $x\in M$.

\noindent So the result of Part 2 of Proposition \ref{ProI} can be interpreted in the following way:

 The distribution $\hat{\cal D}=l^1(\hat{\cal X})$ is minimal among all the $l^1$- characteristic distribution $l^1({\cal Y})$, generated by the family of vector fields ${\cal Y}$ which  satisfies (Hi) and (Hii),  contains ${\cal X}$ and which are ${\cal X}$-invariant.

\end{rem}
\noindent \begin{proof}\so{\it\ Proof}\\
\textbullet$\quad${\bf Proof of part {\it 1}.}\\
 We want to prove that
 $T_z\Phi[{\hat{\cal D}}_z] = {\hat{\cal
D}}_{\Phi(z)}$ for any $z\in$Dom$(\Phi)$ and for any flow $\Phi$ of vector field of $\cal X$ and $\hat{\cal X}$\\

We first show that this is true for any flow $\phi^X_t$ where $X\in {\cal X}$.
Take any $Z\in\hat{\cal X}$ such that $z$ belongs to Dom$(Z)$ and set $x=\phi^X_t(z)$.
There exists $Y\in {\cal X}$ and a finite composition $\Phi$ of flows of vector fields of $\cal X$ such that $Z=\Phi_*(\n Y)$ for some $\n>0$.  So we have $Z'=(\phi^X_t)_*(Z)=\Phi'_*(\n X)$ where $\Phi'=(\phi^X_t\circ\Phi )$. But, there exists  $\n'>0$ such that  $\n' Z$ belongs to $\hat{\cal X}_{x}$, in particular $Z'(x)$ belongs to $\hat{\cal D}_x$. As $\hat{\cal D}_x$ is generated by $\{Y(x),\; Y\in\hat{\cal X}_x\}$
  we then have :
\begin{eqnarray}\label{inclusion}
 T_z\phi^X_t[{\hat{\cal D}}_z] \subset {\hat{\cal D}}_{x}.
\end{eqnarray}
As $(\phi^X_t)^{-1}=\phi^X_{-t}$, by the same argument we get $T_x(\phi^X_t)^{-1}[{\hat{\cal D}}_x] \subset {\hat{\cal D}}_{z}$
and from (\ref{inclusion}) we get
$$T_z \phi^X_t[T_x(\phi^X_t)^{-1}[{\hat{\cal D}}_x]]=\hat{\cal D}_{\phi^X_t(z)}=\hat{\cal D}_x\subset T_z\phi^X_t[{\hat{\cal D}}_z] $$

Now, from (\ref{Yflow}) and the previous argument, we also have $T_z\phi^Y_t[{\hat{\cal D}}_z] = {\hat{\cal D}}_{\phi^Y_t(z)}$ for any $z\in$Dom$(\phi^Y_t)$ and  for any flow $\phi^Y_t$ with $Y\in \hat{\cal X}$.\\


\noindent \textbullet$\quad${\bf Proof of part  {\it 2}.}\\
Let be $x\in M$ and  $Z\in \hat{\cal X}$
such that  $x\in $Dom$(Z)$. As before, we have $Z=\Phi_*(\n Y)$ for some finite composition of flows of vector fields of $\cal X$ and $Y$ is a vector field of $\cal X$ and $\n>0$.
 \begin{displaymath}
\begin{tabular}{ r c l }
$Z(x)$ & $=$& $Z(\Phi(\Phi^{-1}(x)))$\\
&&\\
&$=$& $T_{\Phi^{-1}(x)}\Phi(\n Y(\Phi^{-1}(x)))$
\end{tabular}
\end{displaymath}
As $\D$ is ${\cal X}$-invariant  we obtain that  $Z(x)$ belongs to $\D_x$  and we get $\hat{\cal
D}_x\subset {\D}_x$. In particular, if $\D={\cal D}$,  it is obvious that $\hat{\cal
D}_x$ contains $ {\cal D}_x$, so  we get an equality.\\
This ends the proof of part {\it 2}.\\

\noindent \textbullet$\quad${\bf Proof of part  {\it 3}.}\\
In this proof  we will use some notations and results proved in section \ref{appendice}. In each case, we will mention the  precise references of these notations and results.\\
Let be $x\in M$ for which all assumptions in part {\it 3} are satisfied.
Denote by  $(V_x,\phi)$ the  chart centered at  $x$ such that $\{X_\a,\;\a\in A\}\subset\hat{\cal X}_{V_x}$ Then,  $V_x\subset $Dom$(X_\a)$  for each $\a\in A$ and we set
$$k=\sup \{||J^{s+2}(X_\a){(y)}||_{\phi},\; \a\in A ,
\; y\in V_x\}.$$
Without loss of generality, we can assume that $V_x$ is an open subset $V$ of the Banach space $E\equiv T_xM$ and also that $TM\equiv V\times E$ on $V_x$.  We choose  $r>0$ such that   $B(x,2r)$ is contained in $V$. For the sake of simplicity, we denote
by  $$\{\epsilon_\a=X_\a(x)\}_{\a\in A}$$  the symmetric  unconditional  basis of $\hat{\cal D}_x$

There exists an isomorphism $T$ from $\hat{\cal D}_x$
 to $\R^A$ such that : $T(\epsilon_\a)=e_\a$ where  $\{e_\a\}_{\a
\in A}$  is the canonical basis of $\R^A$. So we can choose  $\r>0$ such that
the image by $T$ of the open  ball $ B(0,\r)\subset \hat{\cal D}_x$ is contained in  $B(0,\frac{r}{k})\subset\R^A$.

 Given any fixed $w=\displaystyle \sum_{\a\in A}t_\a\epsilon_\a
\in B(0,\r)$, we set  $T(w)=\tau=(\t_\a)_{\a \in \a}
$. Of course, $T(w)\in B(0,\frac{r}{k})$.
  By application of
Theorem  \ref{I} on $V$ in the particular case where :

  $\xi= \{X_\a\}_{\a\in A}$
,$\quad I=\R$
, $\quad u =\Gamma^{\tau}$ (see section \ref{appendice} subsection \ref{gamma1}),
 $\quad t_0 =0$
,$\quad T'$ is any real number, large enough,
and $ T_0 =|| \tau||_1$. \\

We have already proved  that :
$$ \Gamma^{\tau} \in L^1_b(\R), \quad \textrm{ with } \quad
||\Gamma^{\tau}||_{\infty}=1.$$
 Let be  $I_0 = [-T_0,T_0]$ and  $B_0 = B(x,
r-kT_0)$. As $T_0 < \displaystyle
\frac{r}{k}$,  there exists a flow $\Phi_{\Gamma^{\tau}}$ defined on $J_0 \times B_0$. From Theorem
\ref{II}, $\Theta= \Psi^x\circ T$, is a map of class
 $C^s$ from  $B(0,\r)\subset \hat{\cal D}_x$
 with values in an open set of $E$ contained in $V$.
We then have:
\begin{eqnarray}\label{defTheta}
\Theta(w)=\Psi^x(\t)=\phi^\xi_\t(x)=\Phi_{\Gamma^{\tau}}(\|\tau\|_1,x)
\end{eqnarray}

The exact  expression of   $\psi^x$ is given
in section \ref{appendice}

It follows from Theorem 2 that $\Theta$ is a map of class $C^s$ with $s>0$ from $B(0,\r)$ into $V$. We can consider $D\Theta_w$ as a field on $B(0,\r)$ of operators from $\hat{\cal D}_x\equiv\R^A$ into $T_xM\equiv E$. On the other hand, we have:

\begin{displaymath}
\begin{tabular}{ r c l }
$D\Theta_{0}(\epsilon_\a)$ & $=$& $D\Psi^x_{(0)}(T(\epsilon_\a))$\\
&&\\
&$=$& $D\Psi^x_{(0)}({e_\a})$\\
&&\\
&$=$& $ \epsilon_\a$
\end{tabular}
\end{displaymath}

\noindent So $D\Theta_0$  is an injective operator from $\hat{\cal D}_x$ into $T_xM$.

Now from \cite{Pe} we have
\begin{lem}\label{inject}${}$
\begin{enumerate}
\item Consider  two Banach spaces ${E_1}$ and  $E_2$ and $i:E_1\ap E_2$ an injective continuous  operator.  Let $\Theta_y$ be a  continuous field of continuous operators  of $L(E_1,E_2)$ on an open neighbourhood $V$ of $x\in E_1$ such that $\Theta_x=i$. Then there exists  a neighbourhood  $W$ in $V$ such that $\Theta_y$ is an injective operator  on $W$.
\item Let $f:U \ap V$ be a  map of class $C^1$ from two open sets $U$ and $V$ in Banach spaces $E_1$ and $E_2$ respectively such that $T_uf$ is injective at $u\in U$.  Then there exists an open neighbourhood $W$ of $u$ in $U$ such that the restriction of $f$ to $W$ is injective.
\end{enumerate}
\end{lem}

By applying this lemma, we conclude that, for $\r$ small enough,  $\Theta: B(0,\r)\ap V$ is a weak submanifold of class $C^s$.\\

It remains to show that $D\Theta_w(\hat{\cal D}_x)= \hat{\cal D}_{\Theta(w)}$.\\

 Given $v=\dis\sum_{\a\in A}v_\a\epsilon _\a$, we set $\s=T(v)$. From (\ref{decompositiondif}) (section \ref{appendice}), we have
$$ D\Theta_w(v)=D\Psi^x_{T(w)}=\D\Psi^x(\tau)\circ{\cal R}(\t)(\dis\sum_{\a\in A}\s_\a \epsilon _\a)$$

On one hand, the map  ${\cal R}(\t)(\dis\sum_{\a\in A}\s_\a \epsilon _\a)=\dis\sum_{\a\in A}\s_\a\D\hat{\Psi}^x_\a((-\tau)^\a)[X_i(\Psi^x_\a(\tau^\a)]$ is a continuous  field  $\t \mapsto {\cal R}(\t)$ of endomorphisms of $\hat{\cal D}_x$ (see Lemma \ref{contR}). As at $\t=0$, the operator  ${\cal R}(0)$ is the identity of $\hat{\cal  D}_x$,  for $\r$ small enough, $w \mapsto {\cal R}\circ T(w)$ is a field of isomorphisms of $\hat{\cal  D}_x$.

On the other hand we have $\Theta(w)=\phi^\xi_{T(w)}(x)$. As $\phi^\xi_{T(w)}$ belongs to ${\cal G}_{\cal X}$, from part {\it 1} of this Proposition,  we have: $D\phi^\xi_{T(w)}(\hat{\cal D}_x)=\hat{\cal D}_{\phi^\xi_{T(w)}(x)}=\hat{\cal D}_{\Theta(w)}$.\\

So we obtain the result required for $\r$ small enough. This ends the proof of part {\it 3}.\\

\noindent \textbullet$\quad${\bf Proof of part  {\it 4}.}\\
The point $x\in M$ for which the assumptions of part {\it 3} of the proposition is true will be  fixed, and    we suppose that $TM$ is trivializable  on the chart domain $V$ (around $x$). We then have:

 \begin{lem}\label{lowtrivial}${}$\\
 Let  $\{X_\a\}_{\a\in A}$ be a family of vector fields on $U\subset V$ which satisfies the condition (LB(s+2)) on $U$ and which is an unconditional symmetric basis of $\hat{\cal D}_x$.
 \begin{enumerate}
\item  There exists a morphism $\Psi:U\times \hat{\cal D}_x\ap TM$ which is a lower trivialization at $x$ such that $\Psi_y(X_\a(x))=X_\a(y)$ for any $\a\in A$.
\item For any integral manifold $f:N\ap U$ of $\hat{\cal D}$ of class $C^s$ through $x$, there exists  a family
$\{Y_\a\}_{\a\in A}$ of vector fields on $N$ defined on a neighbourhood of $f^{-1}(x)$  such $f_*Y_\a=X_\a$ and $\eta=\{Y_\a\}_{\a\in A}$ satisfies the condition (LB(s+2)) at $f^{-1}(x)$.
 \end{enumerate}
 \end{lem}

\noindent \begin{proof}\so{\it Proof}${}$\\
Consider $\tilde{\Psi}:\hat{\cal D}_x\times U\ap \hat{\cal D}$ defined in the following way:

if $w=\dis\sum_{\a\in A}w_\a \epsilon_\a$ we set $\Psi(w,y)=\dis\sum_{\a\in A}w_\a X_\a(y)$.\\

As usual, we set $\tilde{\Psi}_y=\tilde{\Psi}(\;,y)$. Denote by $\bar{\cal D}_y$ the normed subspace defined by $\hat{\cal D}_y$ from the structure of Banach space on $T_yM$, and $i_y: \hat{\cal D}_y\ap T_yM$ the natural inclusion.

At first, as by definition, $\dis\sum_{\a\in A}w_\a$ is absolutely summable,  from the property $LB(s+2)$, it follows that $\Psi(w,y)\in T_yM$ is well defined and $\Psi_y$ is a continuous operator from $\hat{\cal D}_x$ to $\hat{\cal D}_y$ such that $||\Psi_y||\leq K$.  We set $\Psi_y=i_y\circ \tilde{\Psi}_y$. It is clear that  the field $y\ap i_y\circ \tilde{\Psi}(y)$ is smooth.
From this construction, it is easy to see that $\Psi(w,y)=i_y\circ\tilde{\Psi}_y(w)$ is a lower trivialization at $x$ such that $\Psi_y(X_\a(x))=X_\a(y)$ for any $\a\in A$.\\

Let  $f:N\ap U$ be an integral manifold of $\hat{\cal D}$   through $x$ of class $C^s$. Then,  $N$ is a Banach manifold  modeled on the Banach space $\hat{\cal D}_x$. For any open neighborhood $W$ of $x$ the set  $\tilde{W}=f^{-1}(W)$ is an open neighborhood of $\tilde{x}=f^{-1}(x)$. Without loss of generality,  we may assume that $N$ is an open set in $\hat{\cal D}_x$, with $\tilde{x}=0$,  and $M$ is an open set in $E\equiv T_xM$. Modulo these identifications, $f$ is the natural inclusion  of $N$ in $M$, that is the restriction to $N$ of the natural inclusion $i_x:\hat{\cal D}_x\ap T_xM$. In this context, on $i(N)\subset M$, $y\ap \Psi_y$ is a $C^s$  field of continuous linear operators $\textrm{ from }\; \hat{\cal D}_x  \textrm{ into }\;  i_y({\cal D}_y)\equiv i_x(\hat{\cal D}_x)\times\{y\}\subset E\times\{y\}\equiv T_yM$.
From Lemma  2.10 in  \cite{Pe} $y \mapsto \tilde{\Psi}_y$ is also a $C^s$ field of linear operators
$\textrm{ from } \hat{\cal D}_x  \textrm{ into } \hat{D}_y\times\{y\}\equiv i_x(\hat{\cal D}_x)\times\{y\}\equiv T_{y}N$.
It follows that, for any $\a\in A$,  $Y_\a(y)=\tilde{\Psi}(\epsilon_\a)$ is a $C^s$ vector field on $N$ such that $(i_x)_*Y_\a\equiv f_*Y_\a=X_\a$. From the previous  definition of $Y_\a$, it follows that $\eta=\{Y_\a\}_{\a\in A}$ satisfies the condition (LB(s+2)).
\end{proof}

Now we come back to the proof of part {\it 3}. Consider an integral manifold $f:N\ap M$ of $\hat{\cal D}$ of class $C^s$ through $x$ and suppose that the assumption of part {\it 2} is satisfied at $x$.   On $N$, the family  $\eta=\{Y_\a\}_{\a\in A}$ satisfies the condition (LB(s+2)) at $\tilde{x}=f^{-1}(x)$. So for $\r$ small enough,  given any $\t\in B(0,\r)\subset \R^A$,  we can apply Theorem \ref{I} to the family $\eta$ and  $u=\G^\t$ and Theorem \ref{II} on $N$. We then get:

\begin{enumerate}
\item[$\bullet$]a $C^s$ flow $\tilde{\Phi}_{\G^\t}(t,\;)$ of $\tilde{Z}=\dis\sum_{\a\in A}\G_\a^\t Y_\a$ (see section \ref{appendice}) such that  for any $z$ in a small neighborhood $W$ of $\tilde{x}$ we have
$${\Phi}_{\G^\t}(t,f(z))=f\circ \tilde{\Phi}_{\G^\t}(t,z)$$
where  ${\Phi}_{\G^\t}(t,\;)$ is the flow of ${Z}=\dis\sum_{\a\in A}\G_\a^\t X_\a$.
\item[$\bullet$] on $N$, the associated flow   $\tilde{\phi}^\eta_\t({x})=\tilde{\Phi}_{\G^\t}(T,\tilde{x})$
and as in (\ref{defTheta}),  the associated  map\\ $\tilde{\Theta}(w)=\tilde{\phi}^\eta_{T(w)}(\tilde{x})$.\\
\end{enumerate}
  Moreover,  as $T_0\tilde{\Theta}$ is an isomorphism, so for  $\r$ small enough, $\tilde{\Theta}$ is a diffeomorphism from $B(0,\r)$ on a open neighborhood $W$ of $\tilde{x}$ in $N$. On the other hand, from the previous construction, for $\r$ small enough, we have $\Theta=f \circ\tilde{\Theta} $. It follows that $f^{-1}(\Theta(B(0,\r))=\tilde{\Theta}(B(0,\r))=W$. This ends the proof of part {\it 4}.\\
   \end{proof}

 \section{ Structure of weak submanifold on $\cal X$-orbits}\label{weakmani}
  In this section, we will give sufficient conditions under which each $\cal X$-orbits has a structure of weak submanifold of $M$. The first one imposes some  local conditions on the set $\hat{\cal X}$ which leads to integrability of $\hat{\cal D}$ (Theorem \ref{III}) and can  be seen as a generalization of Sussmann's arguments used in \cite{Su}.
The second one   imposes that $\hat{\cal D}$ is  upper trivial and also some local involutivity conditions on  $\hat{\cal X}$.
  \subsection{Structure of manifold and  $\cal X$-orbits}\label{manifolsXorbit}
   Now we will prove some results about integrable  distributions which contain ${\cal D} $ and ${\cal X}$-orbits. This result will be used in each two following subsections.\\

   Consider any set $\cal Y$ of local vector fields which contains $\hat{\cal X}$ and satisfies conditions (Hi). Assume that there exists a weak distribution generated by $\cal Y$:  for instance if $\cal Y$ satisfies (Hii) then we can choose $\D=l^1(\cal Y)$  (see subsection \ref{characdist}).  Assume that $\D$ is integrable on $M$ and for each $x\in M$
there exists a lower trivialization $\Theta :F\times V\ap TM$  for some Banach space $F$ (which depends of $x$) and some neighborhood $V$ of $x$ in $M$.  Let $N$ be the union of all integral manifolds $i_L:L\ap M$ through $x_0$ . Then $i_N:N\ap M$ is the maximal integral manifold of $\D$ through $x_0$ (see Lemma 2.14 \cite{Pe}).

{\it For the clarity of the proof of results in this subsection, for any point $z\in  N$, when $N$ is equipped with the induced topology of $M$, we denote by $\tilde{z}$ the same point of $N$ but when  N is equipped  of its  Banach manifold structure.}

  \begin{prop}\label{varXS}${}$\\
  As previously, let $f\equiv i_N:N\ap M$ be the maximal integral manifold of $\D$  through $x$.
  \begin{enumerate}
 \item Let  $Z\in {\cal X}(M)$ be such that Dom$(Z)\cap  f(N)\not=\emptyset$  and  $Z$ is tangent to $\D$. Set $\tilde{V}_Z=f^{-1}(\textrm{Dom}(Z)\cap  f(N))$. Then $\tilde{V}_Z$ is an open set in $N$ and there exists a vector field $\tilde{Z}$ on  $N$ such that Dom$(\tilde{Z})=\tilde{V}_Z$ and $ f_*\tilde{Z}=Z\circ f $.

Moreover,  if   $]a_x,b_x[$ is the maximal interval on which the integral curve $\g: t\ap \phi^Z(t,{x})$ is defined  in $M$, then the integral curve $\tilde{\g}:t\ap \phi^{\tilde{Z}}(t,\tilde{x})$ is also defined on  $]a_x,b_x[$ and we have
\begin{eqnarray}\label{gZ}
\g=f\circ\tilde{\g}
\end{eqnarray}
\item Let  $\xi=\{X_\b,\; \b\in B\}\subset \hat{\cal X}\subset {\cal Y}$ be which satisfies the conditions (LB(s+2)) on a chart domain $V$ centered at $x\in f(N)$ and consider $\phi^\xi_\t$
for some $\t\in \R^B$ as defined in Theorem \ref{II} and let $\g$ be the $l^1$-curve on $[0,||\t_1||_1]$ associated to $\phi^\xi_\t$ as in Remark \ref{inepB}. Then  there exists a $l^1$-curve $\tilde{\g}:[0,||\t||_1[\ap N$ such that
\begin{eqnarray}\label{l1N}
f\circ \tilde{\g}=\g \textrm{ on } [0,||\t||_1[
\end{eqnarray}
When $\D$ is a closed distribution,  we extend $\tilde{\g}$ to $[0,||\t||_1]$ so that (\ref{l1N}) is true on $[0,||\t||_1]$.
 Moreover, under this last assumption,  to the local diffeomorphism $[\phi^\xi_\t]^{-1}$,  consider the associated $l^1$ curve $\hat{\g}$. Then the curve $\tilde{c}(s)=\tilde{\g}(T-s)$ is a $l^1$-curve  which satisfies (\ref{l1N}) relatively to  $\hat{\g}$.
\end{enumerate}
\end{prop}

\noindent\begin{proof}\so{\it Proof of Proposition \ref{varXS}}\\

$\bullet$ {\bf Proof  of part {\it 1}}

 Fix some $Z\in {\cal X}(M)$ as in Lemma. As $f$ (resp. $T_{\tilde{y}}f$ for any $\tilde{y}\in N$) is injective, there exists a field $\tilde{Z}:\tilde{y}\ap \tilde{Z}(\tilde{y})\in T_{\tilde{y}}N$  such that
\begin{eqnarray}\label{tildeZ}
T_{\tilde{y}}f[\tilde{Z}(\tilde{y})]=Z(f(\tilde{y})), \textrm{ for any } \tilde{y}\in \tilde{V}_Z=f^{-1}[\textrm{Dom}(Z)\cap f(N)]
\end{eqnarray}

 {\it It remains to show  that   the vector field $\tilde{Z}$ is smooth on $\tilde{V}_Z$}.\\

  In fact, it is sufficient to prove this property on some neighborhood  $\tilde{V}$ of any point $\tilde{x}\in \tilde{V}_Z$.  \\
 Note at first that from our assumption about the lower trivialization, we have $\tilde{\D}_x=T_{\tilde{x}}N\equiv F$. So $F$ is independent of $x\in f(N)$.  For any $x\in f(N)$ and an associated  lower trivialization    $\Theta:\R^A\times  V \ap TM$ we will always  choose $V$ such that $TM_{|V}\equiv E\times V$.
 Of course, $f^{-1}(V)$ is an open neighborhood of $\tilde{x}$ in $N$. We also always choose  an open neighborhood $\tilde{V}$ of $\tilde{x}$ in  $f^{-1}(V)$ such that $TN_{| \tilde{V}}\equiv F\times \tilde{V}$.\\

 {\it We assert that   the vector field $\tilde{Z}$ is smooth on $\tilde{V}$}.\\

 Indeed,  from convenient analysis (see \cite{KrMi}), recall that for a map $g$ from an open set $U$ in a Banach space $E_1$ to a Banach space $E_2$ we have the equivalent following properties:
 \begin{enumerate}
 \item[(i)] $g$ is smooth;
 \item [(ii)] for any smooth curve $c:\R\ap U$ the map  $t\mapsto g\circ c(t)$ is smooth;
 \item [(iii)] the map $t\mapsto <\a, g\circ c(t)>$ is smooth for any $\a\in E^*_2$.
  \end{enumerate}
Fix some $\tilde{y}\in \tilde{V}_Z$. As we have already seen,  we can choose a neighborhood $\tilde{V}$ of $\tilde{y}\in \tilde{V}_Z$ such that $TN_{| \tilde{V}}\equiv F\times \tilde{V}$. So, without loss of generality, we can suppose that $\tilde{V}$ is  an open set in $F$ and   $V$ an open set in $E$ and $f\equiv T_{\tilde{y}}f$ on $\tilde{V}$. For simplicity, let be $\theta= T_{\tilde{y}}f: T_{\tilde{y}}N\equiv F\ap T_yM\equiv E$ where $y=f(\tilde{y})$ with our conventions. In these conditions, $\tilde{Z}$ is a map from $\tilde{V}$ to $F$ and $Z$ is  a smooth map from $V$ to $E$. Note that, according to (\ref{tildeZ}), we have
$$\theta\circ\tilde{Z}(\tilde{y})=Z\circ\theta(\tilde{y})$$ for any $\tilde{y}\in \tilde{V}$.
Choose any $\o\in E^*$.  For any smooth curve $c:\R\ap \tilde{V}$, we then have:
$$<\o,Z\circ\theta \circ c>=<\o, \theta\circ\tilde{Z}\circ  c>=<\theta^t(\o), \tilde{Z}\circ  c>$$

As the adjoint  $\theta^t$ of $\theta$ is surjective,  according to the previous argument of convenient analysis we conclude that $\tilde{Z}$ is smooth on $\tilde{V}$.\\

Now, if $x=f(\tilde{x})$, from the relation $ f_*\tilde{Z}=Z\circ f $ we get :
$$\phi^Z(t,x)=f\circ\phi^{\tilde{Z}}(t,\tilde{x})$$

for any $t$ for which $\phi^{\tilde{Z}}(t,\tilde{x})$ is defined. In particular, this relation exists for some interval $]-\varepsilon,\varepsilon[$.  \\

Given the maximal interval $]a_x,b_x[$ as in the Lemma, choose any $\t\in [0,b_x[$. For each $t\in [0,\t]$ we have an integral manifold $f_t:L_t\ap M$ which is an integral manifold of $\D$ through $\phi^Z(t,x)$. As $\phi^{\tilde{Z}}(t+s,x)=\phi^{\tilde{Z}}(t,\phi^{\tilde{Z}}(s,x))$, by the previous argument, there exists some sub-interval on which the curve $s\ap \phi^{\tilde{Z}}(s,x)$ belongs to $L_t$. If we set $L_\t=\dis\cup_{t\in [0,\t]}L_t$, by connexity argument, using Lemma 2.14 \cite{Pe}, it follows that $i_{L_\t}:L_\t\ap M$ is an integral manifold of $\D$ through $x$. But by construction $L_\t$ is an open submanifold of $N$.  It follows that  (\ref{gZ}) is true on $[0,b_x[$; the same arguments works for any $\t\in ]a_x,0]$. This ends the proof of part 1.\\

$\bullet$ {\bf Proof of  part {\it 2}}

Now, let be some $\xi=\{X_\b,\; \b\in B\}\subset \hat{\cal X}\subset {\cal Y}$ satisfying the required conditions. According to Theorem \ref{II}, we have a map $\Psi^x$ from some neighborhood $U$ of $0\in \R^B$ into $V$ of class $C^s$. From part 1, on $N$, we have a family of smooth vector fields $\tilde{\xi}=\{\tilde{X}_\b,\;\b\in B\}$ such that Dom$(\tilde{X}_\b)=\tilde{V}=f^{-1}(V)$ for any $\b\in B$.  Fix some $\t\in U$. According to Remark \ref{inepB},   and (\ref{tildeZ}), by induction, we can construct a curve $\tilde{\g}_\t:[0,||\t||_1[$ such that
\begin{eqnarray}\label{tildeX}
f\circ \tilde{\g}_\t=\Phi^\xi_\t(t,x) \textrm{ for any } t\in [0,||\t||_1[
\end{eqnarray}
Suppose that $\D$ is closed. So $\D_z$ is closed in $T_zM$ for any $z\in N$ and it follows that  the topology of $N$ as weak manifold is nothing but the induced topology of $M$ on $N$. 
The endpoint  $y=\g(||\t||_1)$ belongs to $V$. So $\hat{\g}:[0,||\t||_1]\ap M$ defined by $\hat{\g}(s)=\g(||\t||_1-s)$ is an integral curve of the vector field
$$Z=\dis\sum_{\b\in B}u_\b X_\b$$
 where $(u_\b)=\hat{\G}^\t$ is  associated to $\hat{\phi}^\xi_\t$.\\
 On the other hand, we have an integral manifold $i_L:L\ap M$ through $y$. We choose  a neighborhood $U\subset V$ of $y$ such that we have $TM_{| U}\equiv U\times T_yM$. From our assumption, again,  the topology of $L$ as weak manifold is nothing but the induced topology of $M$ on $L$.

 From  part 1,   $\tilde{U}=U\cap L=(i_L)^{-1}(U\cap L)$  is an open neighborhood of $\tilde{y}=(i_L)^{-1}(y)$ in $L$ and  we have a family $\xi=\{\tilde{Y}_\b,\; \b\in B\}$ such that $(i_L)_*(\tilde{Y}_\b)=[X_\b]_{| \tilde{U}}$.

  From our notations we have $TM_{| U}\equiv U\times T_yM$ and $T_yL$ is a Banach subspace of $T_yM$. So for each $z\in \tilde{U}$ we have an induced  norm on the  finite order jets of vector fields induced from $||.||_\phi$ on the finite jets of vector fields on $\tilde{U}$. As $\{X_\b,\;\b\in B\}$ satisfies the conditions (LB(s+2)) on $V$, and $U\subset V$,  the family $\{\tilde{Y}_\b;\,\b\in B\}$ will  also satisfies  the condition (LB(s+2)) on $\tilde{U}$. So by application of Theorem \ref{I} on $\tilde{U}$ to $Z$ and the unicity of the integral curve through $y$ we have  obtained that $\g(||\t||_1-s)=\hat{\g}(s)$ belongs to $L$ for $0\leq s<\varepsilon$ with $\varepsilon>0$ small enough. We then have  $N\cap L\not=\emptyset$.  It follows that $\tilde {U}$ is an open set of $N$ and in particular  $y$ belongs to $N$ and we can extend $\tilde{\g}$ to $[0,||\t||_1]$

  For the last part, the $l^1$ curve associated to $[\phi^\xi_\t]^{-1}$ is  $\hat{\g}(s)=\g(||\t||_1-s)$ on $[0,||\t||_1]$ and we trivially obtain  the result from the previous proof .\\
\end{proof}

  \subsection{ Structure of weak submanifold on $\cal X$-orbits under local regularity conditions}\label{loccond}

Now we suppose  that $\cal X$ is a set of vector fields on $M$ which satisfies the assumptions (H)=(Hi,Hii,Hiii) that is to say previous conditions  (Hi) and (Hii)  and also the assumption  of Proposition \ref{ProI}, part {\it 3}:\\

(Hiii) there exists a  finite,  countable or  eventually uncountable set $A$ of indexes such that $\hat{\cal D}_x$ is isomorphic to $\R^A$ and a family $\{X_\a, \a\in A\}\subset \hat{\cal X}_x$
such that $\{X_\a, \a\in A\}$ satisfies the condition (LB(s+2)), for some $s>0$,  and $\{X_\a(x), \a\in A\}$ is a symmetric unconditional basis of $\hat{\cal D}_x\equiv \R^A$.\\
\smallskip
\begin{prop}\label{propS} ${}$
\begin{enumerate}
\item For all  $x$ in $M$, the $\hat{\cal
X}$-orbit and the  $\cal
X$-orbit   passing through  $x$ are equal.
 \item The distribution $\hat{\cal D}$ is lower trivial on $M$.
 \item The distribution $\hat{\cal D}$ is integrable. Each  maximal integral manifold of  $\hat{\cal D}$ has a natural  smooth structure of weak connected Banach submanifold, modeled on some $\R^A$
 where $A$ is a finite, countable or  eventually uncountable set of indexes. Moreover, any maximal  integral manifold of $\hat{\cal D}$  is contained in a $\cal X$-orbit.\\
 \end{enumerate}
 \end{prop}

 \begin{theor}\label{III}${}$\\
  If ${\cal  X }$ satisfies the assumptions (H) at each point of $M$,  then $\hat{\cal  D}$  is integrable. Moreover, we have the following properties:
  \begin{enumerate}
  \item[(i)] Each $\cal X$-orbit $\cal O$ is the union of the maximal integral manifolds which meet  $\cal O$ and such an integral manifold is dense in $\cal O$.
  \item[(ii)] Let $\bar{\cal D}$ be  the closed distribution  generated by $\hat{\cal X}$. If $\bar{\cal D}$  is lower trivial and integrable, then, the  $\cal X$-orbit of $x$  is a dense subset in the maximal integral manifold through $x$.
  \item[(iii)]  If $\hat{\cal D}$ is a closed distribution then each ${\cal X}$-orbit is a maximal integral manifold of $\hat{\cal D}$  modeled on some $\R^A$.
  \end{enumerate}
 \end{theor}
 \begin{rem}\label{dimfin}${}$\\
 At any point  $x\in M$ where $\hat{\cal D}$ is a finite dimensional vector space,  $\hat{\cal D}_x$ is isomorphic to some $\R^n$ and  we can always  choose a  finite set $\{X_1,\cdots X_n\}\subset \hat{\cal X}$ such that $\{X_1(x),\cdots X_n(x)\}$ is a basis of $\hat{\cal D}_x$. Moreover for finite set  $\{X_1,\cdots, X_n\}$ we can always find an open neighborhood of $x$  so that  the condition (LB(s+2)) is satisfied on $V$ by this set.
  So, in this case, the assumption (H) is satisfied at $x$. So, if $\hat{\cal D}$ is finite dimensional,  from Theorem \ref{III}  any ${\cal X}$-orbit is a finite dimensional submanifold of $M$.\\
  \end{rem}

\noindent\begin{proof}\so{\it Proof of Theorem \ref{III}}\\
The integrability of $\hat{\cal D}$ is a direct consequence of part 3 of Proposition \ref{propS}.

Moreover, again from part 3 of Proposition \ref{propS} we know that each maximal integral manifold $N$ is contained in a $\cal X$-orbit $\cal O$. It remains to show that such an integral manifold is dense in   $\cal O$. As the binary relation associated to the $\cal X$-orbit is symmetric, $\cal O$ is the $\cal X$-orbit of any point of $\cal O$.  So, if $L$ contains $x$, then, from Proposition \ref{l1Xorbit} part 2 and Proposition \ref{varXS}, any $y\in{\cal O}$ must belong to the closure of $L$ (in $M$). \\

Now, assume that the closed distribution $\bar{\cal D}$ generated by $\hat{\cal X}$ is lower trivial and integrable. Denote by $\cal O$ the ${\cal X}$-orbit of $x$. Choose some $y\in {\cal O}$ and let $\Psi\in {\cal G}_{\cal X}$ be such that $\Psi(x)=y$. According to Comments \ref{Xl1curve}, we can associate to $\Psi$  a finite sequence of points $(x_k)_{k=0,\cdots,n}$ and a finite family $\{\g_k\}_{k=1,\cdots,n}$ of $l^1$-curves associated to some $\phi^{\xi_k}_{\t_k}$ which joins $x_{k-1}$ to $x_k$ and with $x_0=x$ and $x_n=y$. Let $N$ be the maximal integral manifold  of $\bar{\cal D}$ through $x$. As $x_0=x$, according to Proposition \ref{varXS},  there exists  a $l^1$curve $\tilde{\g}_1$  in $N$ such that $i_N\circ\tilde{\g}_1=\g_1$ so $x_1$belongs to $N$. By induction we can construct a $l^1$-curve $\tilde{\g}_k$ in $N$ such that $i_N\circ\tilde{\g}_k=\g_k$ and then $x_k$ belongs to $N$. So, for $k=n$ we obtain that $x_n=y$ belongs to $N$. In particular, by part $1$, each maximal integral manifold $L$ of $\hat{\cal D}$  which meets $\cal O$ is contained in $N$. So, as $\bar{\cal D}$ is closed, the topology of Banach manifold on $N$ is the induced topology as subset of $M$. So, any maximal integral manifold of $\hat{\cal D}$ contained in $\cal O$ is dense in $\cal O$, as subset of $N$. \\
Denote by $\bar{\cal O}$ the closure of $\cal O$ in $N$. So $\bar{\cal O}$ is a connected closed subset of $N$. Consider any $y\in \bar{\cal O}$.
 Let $L$ be the maximal  integral manifold of $\hat{\cal D}$ through $y$. As $L$ is arc-connected and the inclusion of $L$ (with the topology of  Banach manifold ) in $N$ is continuous, it follows that $L$ is contained in $N$. Let ${\cal O}'$ be the $\cal X$-orbit of $y$. From previous arguments, ${\cal O}'$ is also contained in $N$.  Let $(y_k)$ be a sequence in $\cal O$ which converges to $y$. Given any $z\in {\cal O}'$ let be $\Phi\in {\cal G}_{\cal X}$ such that $\Phi(y)=z$.
 Now as each $y_k$ belongs to $\cal O$ and $\cal O$ is invariant by any local diffeomorphism of ${\cal G}_{\cal X}$, it follows that  $z_k=\Phi(y_k)$, for $k$ large enough. So $z=\dis\lim_{k\ap\infty}z_k$, and then $z$ belongs to $\bar{\cal O}$. Finally  we get ${\cal O}'\subset \bar{\cal O}$ and, in  particular, the maximal integral manifold $L'$ of $\hat{\cal D}$ through $y$ is contained in $\bar{\cal O}$.\\
  On the other hand,
as $T_yN$ is the closure in $T_yM$  of the normed subspace $T_yL'$, there exists a neighborhood  $U$ of $y$ in $L'$ (for the two  topologies  on   $L$) such that, the closure of ${U}$ in $N$   is a closed set of $N$ with non empty interior. But $U\subset L'\subset \bar{\cal O}$, so it follows that $\bar{\cal O}$ is open. By connexity argument, we get $\bar{\cal O}=N$.\\

Now if  $\hat{\cal D}$ is a closed distribution, obviously we have $\bar{\cal D}=\hat{\cal D}$ so, the assumptions of property  are satisfied. So part (iii) is a direct consequence of properties  (i) and (ii).
\end{proof}

\noindent\begin{proof}\so{\it Proof of Proposition \ref{propS}}\\


\noindent \textbullet$\quad${\bf Proof of part {\bf 1}}.\\
This result comes from  ${\cal G}_{ \cal X}
={\cal G}_{ \hat {\cal X}}$   (Proposition \ref{l1Xorbit} part {\it 3}).\\

\noindent \textbullet$\quad${\bf Proof of part {\bf 2}}.\\
This result is a consequence of Lemma \ref{lowtrivial} part {\it 1}\\

\noindent \textbullet$\quad${\bf Proof of part 3}.\\

From Proposition \ref{ProI} part {\it 3}, for any $x\in M$ we have a $C^s$ integral manifold through $x$, with $s\geq 1$. As $\hat{\cal D}$ is a lower trivial weak distribution, consider the set
$${\cal X}^{-}_{\cal D}=\{X(u)=\Psi_x(u,y) \textrm{ for any lower  trivialization } \Psi_x:\hat{\cal D}_x\times V\ap TM \textrm{ and any } x\in M\} $$
  As through $x$, we have an integral manifold of class $C^s$, $s\geq 1$, from the proof of Proposition 2.8 in \cite{Pe}, it follows that ${\cal D}$ is ${\cal X}^{-}_{\cal D}$-invariant. So from Theorem 1 of \cite{Pe} we have a smooth integral manifold through $x$.
Moreover, if we consider  the following equivalence relation on $M$:
$$x{\cal R}y \textrm{ iff there exists an integral manifold of } \hat{\cal D} \textrm { passing  through } x \textrm{ and } y$$
then each equivalence class $L$ has a natural structure  of weak Banach submanifold modeled on $\hat{\cal D}_x$ for any $x\in L$ and $L$ is an integral manifold of $\hat{\cal D}$.
Take such an equivalence class $L$   and denote by $i_L$ the natural inclusion of $i_L$ of $L$ (endowed with its Banach structure) into $M$.  From Proposition \ref{ProI} part {\it 3}, for any $x\in L$, there exists an open ball $B(0,\r_x)\subset \R^A\equiv \hat{\cal D}_x$ and a $C^s$ map $\Theta_x: B(0,\r_x)\ap M$ which is a $C^s$ integral manifold of $\hat{\cal D}$ through $x$ and such that $\Theta_x(B(0,\r_x))$ is an open set of $L$.  So, $P_x=\Theta_x(B(0,\r_x))$ has an induced structure of smooth Banach manifold modeled on $\hat{\cal D}_x$ (isomorphic to $\R^A$ for some appropriate set of indexes $A)$. In particular, the natural inclusion $i_x : P_x\ap M$ is a smooth integral manifold of $\hat{\cal D}$ through $x$.
Now take some  $x\in L$. For any $y\in L$ we have a continuous curve $\g:[a,b]\subset \R\ap L$ such that $\g(a)=x$ and $\g(b)=y$. By compactness  of $\g([a,b])$ we have a finite  covering of $\g([a,b])$ by a family of open sets $\{ \Theta_{x_i}(B(0,\r_{x_i})\}_{i=1,\cdots, n}$ such that $x_i\in \g([a,b])$, $x_1=x$ and $x_n=y$. Now choose any $y_i\in \Theta_{x_i}(B(0,\r_{x_i}))\cap\Theta_{x_{i+1}}(B(0,\r_{x_{i+1}}))\cap  \g([a,b]) $ for $i=1,\cdots,n-1$. From the construction of each $\Theta_x$, there exists $\Phi_i\in {\cal G}_{\cal X}$ (resp. $\Phi'_i\in {\cal G}_{\cal X}$) such that  $\Phi_i(x_{i})=y_i$ (resp. $\Phi'_i(x_{i+1})=y_i$). So the composition:
\begin{eqnarray}\label{finitcompo}
\Phi=\Phi_1\circ [\Phi'_1]^{-1}\circ\cdots\circ \Phi_{n-1}\circ [\Phi'_{n-1}]^{-1}
\end{eqnarray}
is an element of ${\cal G}_{\cal X}$ such that $\Phi(x)=y$. It follows that $L$ is contained in the ${\cal X}$-orbit of $x$.

\end{proof}

  \subsection{ Structure of weak submanifold on $\cal X$-orbits under involutivity conditions}\label{stronginvol}

A weak distribution ${\D}$ is called  {\bf (locally) upper  trivial} (upper trivial for short)  if,  for each $x\in M$, there exists an open neighborhood $V$ of $x$, a Banach space ${F}$ and  a smooth map $\Phi:F\times  V \ap TM$  (called {\bf upper trivialization})  such that :
\begin{enumerate}
\item[(i)] for each $y\in V$,  $\Phi_y\equiv \Phi(\;,y):F\ap T_yM$ is a continuous operator with $\Phi_y(F)={\D}_y$;
\item [(ii)] $\ker \Phi_x$  complemented in $F$;
\item[(iii)] if $F= \ker \Phi_x\oplus S$,  the restriction $\theta_y$  of $\Phi_y$ to $S$  is injective for any $y\in V$;
 \item[(iv)]  ${\Theta}(u,y)=({\theta}_y\circ [{\theta}_x]^{-1}(u), y)$ is  a lower trivialization of ${\cal D}$.
 \end{enumerate}
In this case the map  $\Theta$ is called the {\bf associated lower trivialization}.

 In this case, each lower section $X_v=\Theta(v,\;)$  with $v\in\D_x$ can be written as $X_v=\Theta(\Phi(v',x),\;)$ for any $v'\in F$ such that $\Phi(v',x)=v\in \D_x$.

An  upper trivial weak distribution ${\D}$ is called {\bf Lie bracket invariant}  if, for any $x\in M$, there exists an upper trivialization $\Phi:F\times V\ap TM$ such that, for any $u\in F$ , there exists $\varepsilon>0$,
 and,  for all $0<\t<\varepsilon$ there exists  a smooth field of operators $C:[-\t,\t]\ap L(F,F)$ with the following property
 \begin{eqnarray}\label{condLieinvst}
 [X_u,Z_v](\g(t))=\Phi(C(t)[v],\g(t)) \textrm{ for any } Z_v=\Phi(v,\;) \textrm{ and any  } v\in F
\end{eqnarray}
along the integral curve $t\mapsto \phi^{X_u}_t(x)$ on $[-\t,\t]$ of  the lower section $X_u=\Theta(\Phi(u,x),\;)$ .

With these definitions we have:

\begin{theo}\label{lieinv}${}$\\
Let ${\D}$ be an  upper trivial weak  distribution. Then $\D$ is integrable if and only if ${\D}$ is Lie bracket invariant.\\
\end{theo}
\bigskip

  We now  come back to our original context. Consider any set $\cal Y$ of local vector fields which contains $\cal X$ and which satisfies properties (Hi) and (Hii). We have seen that if $\L$ is any ordered set of indexes of same cardinal as the set $${\cal Y}_x=\{Y\in {\cal Y}  \textrm{ such that }   x\in\textrm {Dom}(Y)\}$$ then we have a surjective linear map:
$T:l^1(\L)\ap l^1({\cal Y})_x$.\\  Let $\D$ be the weak distribution $l^1({\cal Y})$ and   index the set ${\cal Y}_x$ as set $\{Y_\l, \L\in \L\}$.
Assume that $\D$ has the following properties labelled (H'):
\begin{enumerate}
\item[(H'1)] for any $x\in M$ there exists an upper trivialization $\Phi:l^1(\L)\times V\ap TM$ such that $\Phi(e_\l)=Y_\l$ for each $\L\in \L$ where $\{e_\l\}_{\l\in\L}$ is the canonical basis of $l^1(\L)$;
\item[(H'2)] for any $x\in M$ there exists a neighborhood $V$ of $x$  such that $V\subset \dis\cap_{\l\in \L}\textrm{Dom}(Y_\l)$,  and a constant $C>0$ such that  we have
\begin{equation}\label{hypinvol}
[Y_\l,Y_\m](y)=\dis\sum_{\n\in \L}C_{\l\m}^\n(y) Y_\n(y) \textrm{ for any } \l,\m\in \L
\end{equation}
where each $C_{\l\m}^\n$ is a smooth function  on $V$, for any  ${\l,\m,\n\in \L}$ and  we have $$\dis\sum_{\n\in A}|C_{\l\m}^\n(y)|\leq C$$ for any $y\in V$.
 \end{enumerate}

  \begin{theor}\label{III'}${}$
  \begin{enumerate}
  \item Under the previous assumptions (H'), the distribution $\D$ is integrable.
 \item If $\D$ is an integrable distribution which satisfies assumption (H'1), then  $\hat{\cal D}_x$ is contained in $\D_x$  for any $x\in M$. Moreover if $\D$ is closed then  each  ${\cal X}$-orbit  is contained in a maximal integral manifold of $\D$.
  \end{enumerate}
  \end{theor}

 To the set ${\cal X}$ we can associate the sequence of families
 $${\cal X}={\cal X}^1\subset {\cal X}^2={\cal X}\cup\{[X,Y],\; X,Y\in {\cal X}\}\subset \cdots\subset {\cal X}^k={\cal X}^{k-1}\cup\{[X,Y],\; X\in {\cal X}, Y\in {\cal X}^{k-1}\}\subset\cdots$$
The set  ${\cal X}^k$ always satisfies the condition (Hi). Moreover, if it satisfies condition (Hii),  the distribution ${\cal D}^k= l^1({\cal X}^k)$ is well defined.

\noindent By application of the previous  result to   $\D=\hat{\cal D}$ or $\D={\cal D}^k$ we get:

 \begin{theor}\label{III"}${}$
 \begin{enumerate}
 \item If the distribution $\hat{\cal D}$ satisfies the assumptions (H'), then $\hat{\cal D}$ is integrable and we have the following properties:
 \begin{enumerate}
  \item[(i)] Each $\cal X$-orbit $\cal O$ is the union of the maximal integral manifolds which meet  $\cal O$ and such an integral manifold is dense in $\cal O$.
  \item[(ii)] Assume that the closed distribution $\bar{\cal D}$ generated by $\hat{\cal X}$ is lower trivial and integrable. Then the  $\cal X$-orbit of $x$  is a dense subset in the maximal integral manifold through $x$.
  \item[(iii)]  If $\hat{\cal D}$ is a closed distribution then each ${\cal X}$-orbit is a maximal integral manifold of $\hat{\cal D}$ modeled on some $\R^A$.
  \end{enumerate}
  \item If ${\cal X}$ satisfies  (LBs),  and  if ${\cal D}^k$  satisfies assumptions (H')  for some  $k\leq s$, then  we have ${\cal D}^k=\hat{\cal D}$  and ${\cal D}^k$ is integrable. Moreover, ${\cal D}^k$ satisfies all the previous properties (i), (ii) and (iii).
\end{enumerate}
 \end{theor}

 \begin{exe}\label{serpent}${}$\\
 As in Example \ref{ex} (3), consider a finite  family ${\cal X}=\{X_1,\cdots, X_n\}$ of global vector fields on $M$. We have seen that  the condition (LBs) is satisfied for any $s>0$. Then each set ${\cal X}^k$ is finite and then,  it is clear that each distribution ${\cal D}^k$ is  upper trivial:

 if $n_k$ is the cardinal of  ${\cal X}^k$, we can order ${\cal X}^k$ in a sequence $\{Z_1,\cdots Z_{n_k}\}$ and on each open set $V$ according to the identification $TM\equiv V\times T_xM$ we can consider the upper trivialization $\Phi: V\times\R^{n_k}\ap TM$ defined by:

 $\Phi(y,(t_1,\cdots ,t_{n_k}))=\dis\sum_{i=1}^{n_k} t_i Z_i(y)$; in fact it is an upper trivialization.

 Suppose that the condition (H'2) for ${\cal X}^k$ is satisfied, then from Theorem \ref{III"}, the closed distribution  ${\cal D}^k$ is integrable, and each maximal integral is a Banach submanifold of $M$ which is also a ${\cal X}$-orbit.

 The reader can find such a context in \cite{Ro} where $M$ is the set (denoted   "$\rm{Conf}$") of "configurations"  of the snake (which is a Banach manifold), ${\cal X}$ is the set
 of global vector fields $\{\xi_1,\cdots,\xi_d\}$ on $\rm{Conf}$ (in notations \cite{Ro}). Then ${\cal X}^1$ satisfies the condition (H'2).  Each ${\cal X}$-orbit is nothing but an orbit of the action  \rm{M\" ob} on $\rm{Conf}$ (see \cite{Ro}). From Theorem \ref{III"} we directly obtain  that each orbit is a closed (finite dimensional) submanifold of $\rm{Conf}$.\\

 In \cite{PS}, the reader can find a generalization of the results of \cite{Ro} in the context of Hilbert space and get an application of  the previous result for  a countable set ${\cal X}$ of global vector fields on a Banach manifold.
 \end{exe}

  \noindent \begin{proof}\so{\it Proof Theorem \ref{III'}}\\

  \noindent$\bullet$ {\bf Proof of  part 1}

  According to Theorem \ref{lieinv}, it is sufficient to show that $\D$ is Lie bracket invariant. So  fix some $x\in M$ and consider an upper trivialization $\Phi:l^1(\L)\times V\ap TM$ as in the previous assumption. As $\ker \Phi_x$ is complemented, we have $l^1(\L)=\ker \Phi_x\oplus S$. So  there exists a family $\{\epsilon_\a\}_{\a\in A}$ (resp. $\{\epsilon'_\b\}_{\b\in B}$) of $l^1(\L)$ which is a normalized symmetric unconditional basis of  $S$ (resp. $\ker \Phi_x$) (see Remark \ref{compl1}).
  Now, the canonical unconditional basis $\{e_\l\}_{\l\in \L}$ has a (unique) decomposition:
   \begin{equation}\label{decompocanon}
e_\l=\dis\sum_{\a\in A}f_\l^\a\epsilon_\a+\dis\sum_{\b\in B}{f'}^\b_\l\epsilon'_\b
\end{equation}
such that $\dis\sum_{\a\in A}|f_\l^\a|\leq 1$ and $\dis\sum_{\b\in B}|{f'}^\b_\l|\leq 1$ for any $\l\in \L$.\\

Any lower section can be written as $X_u=\Phi(u, .)$ for some $u=(u_\l)\in  l^1(\L)$.  Such  a section can be written
  $$X_u=\dis\sum_{\l\in \L} u_\l Y_\l$$
 On the other hand consider a neighborhood $V'$ of $x$ in which (H'2) is true and the neighborhood $V\cap V'$ (again denoted by $V$).
 As previously fix some lower section $X_u=\Phi(u, .)$  and
consider $\varepsilon>0$ such that the  integral curve $\g(t)=\phi^{X_u}_t(x)$ is defined on $]-\varepsilon,\varepsilon[$ in $V$.  According to (H'2),  for any $0<|\t|<\varepsilon$,  we define $C:[-\t,\t]\ap L(l^1(\L),l^1(\L))$ in the following way:

 $C(t)[v]=\dis\sum_{\l,\m,\n\in \L} C_{\l\m}^\n(\g(t)) u_\l v_\m e_\n$

 where $v=(v_\m)\in l^1(\L)$.
So from assumption (H'2), we have :

 $||C(t)[v]||\leq C\dis\sum_{\l,\m \in \L}|u_\l||v_\m|\leq C[\dis\sum_{\l\in \L}|u_\l|][\dis\sum_{\m\in \L}|v_\m|]=C||u||_1||v||_1$

 then $C(t)$ is a field of continuous endomorphisms of $l^1(\L)$.\\ On the other hand, for any $v=(v_\m)\in l^1(\L)$, we have

 $Z_v=\Phi(v,.)=\dis\sum_{\m\in \L}v_\m Y_\m$

So we get
$[X_u,Z_v](\g(t))=\Phi(C(t)[v],\g(t))$

\noindent From Theorem \ref{lieinv}  it follows that $\D$ is integrable.\\

 \noindent$\bullet$ {\bf Proof of part  2}${}$\\

    Now suppose that $\D$ is integrable. Fix some $x\in M$  and let  $f\equiv i_N:N\ap M$ be a maximal integral manifold through $x$. We want to show that $\hat{\cal D}_x$ is contained in $\D_x$. It is sufficient to prove that  for any $Y\in \hat{\cal X}_x$, $Y(x)$ belongs to $\D_x$.  For such a vector field there exist vector fields $X_1,\cdots X_p,X\in {\cal X}$ and $\n >0$ such that
  $$Y=(\phi^{X_p}_{t_p}\circ\cdots\circ \phi^{X_1}_{t_1})_*(\n X)$$
   Let $z=(\phi^{X_p}_{t_p}\circ\cdots\circ \phi^{X_1}_{t_1})^{-1}(x)$ be.
   Consider the integral curve $\g_1$ of $X_1$ through $z$: $\g_1(t)=\phi^{X_1}_{t}(z)$ for $t\in [0,t_1]$. As ${\cal X}\subset {\cal Y}$, from Proposition \ref{varXS}, we have a curve $\tilde{\g}_1:[0,t_1]\ap N$ such that $f\circ\tilde{\g}_1=\g_1$, and for any $s\in [0,t_1]$, a neighborhood $\tilde{V}_s$ of $\tilde{\g}_1(s)$ in $N$, and a  vector field $\tilde{Y}_s$ on $\tilde{V}_s$ such that $f_*\tilde{Y}_s=X_1$. In particular we also have
   \begin{eqnarray}\label{compatflow}
f\circ\phi^{\tilde{Y}_s}_r(\tilde{\g}_1(s))=\phi^{X_1}_r(\g_1(s)) \textrm { for any } r \textrm{ small enough}
\end{eqnarray}
Moreover,  we can  find $\tilde{X}$ on the  neighborhood $\tilde{V}_0$ of $\tilde{z}$  such that $f_*(\tilde{X})=\n X$, after having restricted $\tilde{V}_0$ if necessary. By compactness, we can cover $\tilde{\g}_1([0,t_1])$ by a finite number  $\tilde{V}_{s_0},\cdots \tilde{V}_{s_m}$. On $\tilde{V}_{s_0}$ we have $T_{\tilde{z}}[\phi^{\tilde{Y}_{s_0}}_t](\tilde{X}((\tilde{z}))$ which belongs to $T_{\tilde{\g}_1(t)}N=\tilde{D}_{\tilde{\g}_1(t)}$ for any $t$ so that $\tilde{\g}_1(t)$ belongs to $\tilde{V}_{s_0}$. From properties of $\tilde{Y}_{s_0}$ and $\tilde{\g}_1$ , it follows that
\begin{eqnarray}\label{app}
[(\phi^{X_1}_t)_*(\n X)](\g(t))\textrm{  belongs to } \D_{\g(t)}.
\end{eqnarray}
Choose $\s_1$ such that $\tilde{\g}_1(\s_1)$ belongs to $\tilde{V}_{s_1}$. So we have (\ref{app}) for $t=\s_1$. By applying the same argument to $T_{\tilde{z}}[\phi^{\tilde{Y}_{s_0}}_{\s_1}](\tilde{X}((\tilde{z})) $ by choosing $\s_2$ such that $\tilde{\g}_1(\s_2)$ belongs to $\tilde{V}_{s_1}\cap \tilde{V}_{s_2}$, we obtain (\ref{app}) for $t=\s_2$. Finally, by induction we get (\ref{app}) for $t=t_1$. Then by same argument applied to $[(\phi^{X_1}_t)_*(X)](\g(t_1))$ instead of $(\n X)(x)$  and along the curve $\g_2(t_1+t)=\phi^{X_2}_t(\g_1(t_1))$ we obtain that
$$(\phi^{X_2}_{t_2}\circ\phi^{X_1}_{t_1})_*(\n X(x)) \textrm{  belongs to } \D_{\g(t_2)}$$
Again by induction, on $i=2,\cdots, p$, we finally obtain that $Y(x)=(\phi^{X_p}_{t_p}\circ\cdots\circ \phi^{X_1}_{t_1})_*(\n X(z))$ belongs to $\D_x$.

Now we assume that $\D$ is a closed integrable manifold.  Take $x\in M$ and again let be $f=i_N:N\ap M$ the maximal integral manifold trough $x$. We want to show that for any $\Psi\in {\cal G}_{\cal X}$, the point $y=\Psi(x)$ belongs to $f(N)$. From the previous proof we also have obtained that if $\Psi$  is a finite composition $(\phi^{X_p}_{t_p}\circ\cdots\circ \phi^{X_1}_{t_1})$, then $y$ belongs to $N$ . So, from  (\ref{Yflow}), for  $Y\in \hat{\cal X}$ and any $\t\in \R$, $\phi^X_\t(x)$ belongs to $f(N)$ (even when $\D$ is not closed).\\

 Suppose that  $\;\Psi$ is reduced to some $\phi^\xi_\t$, with $\xi=\{Y_\d, \d\in D\}\subset\hat{\cal X}$  and $\t\in \R^D$\\
  Let $\g:[0,||\t||_1]\ap M$ be the curve $\g(t)=\Phi^\xi_\t(t,x)$ where $\Phi^\xi_\t(t,.)$ is the flow associated to $\xi$, $\t$ and $u=\G^\t$ (see Remark \ref{inepB} {\it 1.})
From Proposition \ref{varXS}, part 2  there exists a $l^1$-curve $\tilde{\g}:[0,||\t||_1]\ap N$ such that
\begin{eqnarray}\label{gN}
f\circ \tilde{\g}=\g \textrm{ on } [0,||\t||_1]
\end{eqnarray}
As $\phi^\xi_\t(x)=\Phi^\xi_\t(||\t||_1,x)$, we obtain that $y=\phi^\xi_\t(x)$ belongs to $f(N)$.

For the case $\Psi=[\phi^\xi_\t]^{-1}$,  set again $y=[\phi^\xi_\t]^{-1}(x)$ and let $\g:[0,T]\ap M$ be the $l^1$-curve associated to $\phi^\xi_\t$ and we use the previous notations. Then the $l^1$ curve associated to $\Psi$ is $\hat{\g}(s)=\g(T-s)$ which satisfies $\hat{\g}(0)=x$ and $\g(||\t||_1)$ and $\hat{\g}(||\t||_1)=y=\g(0)$.  From (\ref{gN}), we obtain that $x$ belongs to the maximal  integral manifold through $y$ which, by maximality, is $N$.\\

 In the general case we have $\Psi=\phi_n\circ\cdots\circ\phi_1$  where each $\phi_k$ is a local diffeomorphism of  type $\phi^{Y_k}_{\t_k}$ or $\phi^{\xi_k}_{\t_k}$ or $[\phi^{\xi_k}_{\t_k}]^{-1}$ for $k=1\cdots n$ and all these vector fields belong to $\hat{\cal X}$.
So, by finite induction on $k$, using the previous  partial results, we get the proof of part {\it 2} in the general case.
 \end{proof}\\

\noindent \begin{proof}\so{\it Proof Theorem \ref{III"} }${}$\\

\noindent$\bullet$ {\bf Proof of part 1}\\
By application of Theorem \ref{III'}, part 1 to $\hat{\cal D}$,  it follows that $\hat{\cal D}$ is integrable. We must show that each maximal integral manifold which meets a $\cal X$-orbit $\cal O$ is contained in $\cal O$.

Fix some maximal integral manifold $i_N:N\ap M$ of $\hat{\cal D}$. Fix some $x\in N$ and  consider an upper trivialization $\Phi:\R^\L\times V\ap TM$ as in assumption (H'1). From this assumption, after restricting $V$ if necessary, the set $\xi=\hat{\cal X}_x$ satisfies the condition (LBs)  at any point of $V$ and for any $s\in \N$ (see Example \ref{ex} 2). On the other hand, according to  Lemma 2.10  in \cite{Pe}, we have a neighborhood $\tilde{V}$ of $x$,  for the Banach structure of $N$, so that we have a smooth field of continuous operators $y\ap \tilde{\Phi}_y$ from $\R^\L$ to $T_yN$ such that $\Phi_y(.)=Ti_N\circ \tilde{\Phi}_y$ on $\tilde{V}$.
From Proposition \ref {varXS}, for each $\l\in \L$ we have a smooth vector field on $\tilde{Y}_\l$ such that
 \begin{eqnarray}\label{relY}
Y_\l=(i_N)_*\tilde{Y}_\l \textrm{ on } \tilde{V}
\end{eqnarray}
Note that,  according to the notation used in the proof of Theorem \ref{III'} part 1,  in fact we have $Y_\l(y)=\tilde{\Phi}_y(e_\l)$. So, as previously,  after restriction of $\tilde{V}$ if necessary, the set $\tilde{\xi}=\{\tilde{Y}_\l,\;\l\in \L\}$ satisfies the condition (LB(s+2)) for any $s\in \N$. Applying Theorem\ref{II} to $\tilde{\xi}$ we get a map $\tilde{\Psi}^x: B(0,r)\subset \R^A\ap L$ of class  $C^s$. By the same argument applied to $\xi=\{Y_\l,\;\l\in\L\}$ on $M$, we get a map $\Phi^x:B(0,r')\ap M$ which is of class $C^s$. Using (\ref{relY}) we have $\Psi^x=i_N\circ \tilde{\Phi}^x$ on some $B(0,\r)$ with $\r$ small enough. The linear map  $T_0\tilde{\Phi}^x$ is surjective  and its  kernel is $\ker\Phi_x$. So, for $\r$ small enough, $\tilde{\Psi}^x$ is a submersion and in particular, $\tilde{P}(x,\r)=\tilde{\Psi}^x(B(0,\r))$ is an open set in $N$ (with it Banach structure). If we set  $P(x,\r)=\Psi^x(B(o,\r))$ by definition of a $\cal X$-orbit, the set $P(x,\r)$ is contained in $\cal O$. But, by construction we have $P(x,\r)=i_N(\tilde{P}(x,\r))$ and then we have an open neighborhood $\tilde{P}(x,\r)$ of $x$ (for the Banach structure of $N$) such that $i_N(\tilde{P}(x,\r))\subset {\cal O}$. As we can cover $N$ by such open subsets  and $\cal O$ is the $\cal X$-orbit of any $y \in{\cal  O}$,  we get $N \subset {\cal  O}$. For the density of $N$ in $\cal O$,  we use the same arguments as in the proof of Theorem \ref{III}. The properties (ii) and (iii) have same proofs as in  Theorem \ref{III}.\\


\noindent$\bullet$ {\bf Proof of part 2 }${}$\\
From Theorem \ref{III'}  applied to $\D={\cal D}^k$  we obtain that ${\cal D}^k$ is integrable and, for any $x\in M$  each ${\cal D}^k_x$  contains $\hat{\cal D}_x$. According to part {\it 1} of Theorem \ref{III"}, it remains to show that $\hat{\cal D}_x$ contains ${\cal D}^k_x$  for any $x\in M$.\\
Given $x\in M$, we can suppose that the upper  trivialization  $\Phi: \tilde{\cal D}^k_x\times  V \ap TM$ on a neighborhood $V$ of $x$ is such that $TM_{|V}\equiv E\times V$.  Take any $X\in {\cal X}$ and $Y\in \hat{\cal X}$ so that $x$ belongs to the domain of $X$ and of $Y$. For $0<t<\varepsilon$ small enough so that  the flow $\phi^X_t$ is defined on some neighborhood $U\subset V$ of $x$, we consider the curve $t\ap \dis\frac{1}{t}\{([\phi^X_t]_*Y)_x-Y_x\}$ in $E$. As $\hat{\cal D}$ is ${\cal X}$-invariant, the previous  curve belongs to $\hat{\cal D}_x$, as Banach space. But we have:
$$[X,Y]_x=\dis\lim_{t\ap 0}\dis\frac{1}{t}\{([\phi^X_t]_*Y)_x-Y_x\}$$

As ${\cal D}^k$ satisfies the assumption (H'), the structure of Banach space for ${\cal D}^k_x$ is isomorphic to some $\R^A$. So ${\cal D}^k_x$ has the Schur Property. By using an  argument of weak convergency and Schur's property, $[X,Y]_x$ belongs to $\hat{\cal D}_x$. Now by induction, applying this result for $Y\in {\cal X}^{k-1}$, we obtain the inclusion ${\cal D}^k_x \subset \hat{\cal D}_x$ for any $x\in M$.
\end{proof}

\section{Applications}
\subsection{Criteria of integrability for $l^1$-distribution}
 In this  subsection we will give a criterion  of integrability  for $l^1$-distributions generated by   sets ${\cal X}$  of vector fields on $M$  which satisfies   the assumption (H). We have the following result:

\begin{theor}\label{IV}${}$
\begin{enumerate}
\item Let  $\cal D$ be a  $l^1$-distribution generated by a set of (local) vector fields   $\cal X$  on a Banach manifold $M$ which satisfies the  assumptions (H).  Then  $\cal D$ is lower trivial. Moreover, $\cal D$   is integrable if and only if  $\cal D$ is $\cal X$-invariant.
\item Let  $\cal D$ be a lower trivial $l^1$-distribution on a Banach manifold $M$.  Then there exists generating sets    ${\cal X}$  of ${\cal D}$ which satisfies assumption {\rm{(Hi)}}, { \rm{(Hii)}} and  { \rm{(Hiii)}} . Given any such  generating set $\cal X$ of $\cal D$, then $\cal D$ is integrable if and only if $\cal D$ is $\cal X$-invariant.
\end{enumerate}
\end{theor}

\begin{rem}${}\label{intl1dis}$\\
As any $l^1$-distribution $\cal D$ is a weak distribution, from Theorem 1 of \cite{Pe}, when  $\cal D$ is lower trivial, it is integrable if and only if it is ${\cal X}^-_{\cal D}$-invariant (${\cal X}^-_{\cal D}$ in the set of lower sections of $\cal D$ see \cite{Pe}). 
So,  for lower trivial  $l^1$-distribution, the Theorem 4 gives a necessary and sufficient   condition of integrability for {\bf any}  generating set of $\cal D$ satisfying  {\rm{(Hi)}}, { \rm{(Hii)}} and  { \rm{(Hiii)}}. Note that, if $\cal D$ is finitely generated at each point, these conditions  are automatically satisfied. We then get a generalization of the famous criterion of integrability of Nagano-Sussmann in this context of Banach manifold for finite dimensional distribution. In this  sense, Theorem 4 can be considered as a generalization of this  Nagano-Sussmann's  result in infinite dimension.\\
\end{rem}
 In the Example \ref{ex} 1,  if  the set $\{T(x_\a)\}_{\a\in A}$ is  a family of linearly independent vectors, the conditions of Theorem 4 are satisfied. Of course, this result can be proved directly in an obvious way. Each leaf is the affine space in $E$  associated to the $l^1$  normed space generated by ${\cal X}_0$. On the other hand in the Example \ref{ex} 2, even in analogue conditions, the characteristic distribution of  $\cal X$ is not $\cal X$-invariant. Such a sufficient conditions will be carried by $\Psi$ (see Theorem 4 in \cite {Pe}).\\

\noindent\begin{proof} \so{\it Proof of Theorem \ref{IV}}\\

\noindent {\bf Part 1}

\noindent
From Proposition \ref{ProI} we have ${\cal D}=\hat{\cal D}$ if and only if  $\cal D$ is $\cal X$-invariant.  On the other hand,  $\cal X$ satisfies the assumption (H) (of subsection \ref{weakmani}). By application of Theorem 3, we obtain the first part.\\

\noindent {\bf Part 2}

\noindent
Fix some $x\in M$. From the property of lower triviality, there exists an open neighborhood $V$ of $x$ in $M$,  a smooth map $\Psi:\tilde{\cal D}_x\times  V \ap TM$   such that :
\begin{enumerate}
\item[(i)]  $ \Psi(\tilde{\cal D}_x\times\{y\})\subset {\cal D}_y$ for each $y\in V$

\item[(ii)] for each $y\in V$,  $\Psi_y\equiv \Psi(\;,y):\tilde{\cal D}_x\ap T_yM$ is a continuous operator  and $\Psi_x:\tilde{\cal D}_x\ap T_xM$  is the natural inclusion $i_x$

\item [(iii)] there exists a  continuous operator $\tilde{\Psi}_y: \tilde{\cal D}_x\ap \tilde{\cal D}_y$ such that $i_y\circ \tilde{\Psi}_y=\Psi_y$, $\tilde{\Psi}_y$ is an isomorphism from $\tilde{\cal D}_x$ onto ${\Psi}_y(\tilde{\cal D}_x)$
and  $\tilde{\Psi}_x$ is the identity of $\tilde{\cal D}_x$.
\end{enumerate}

As $\tilde{\cal D}_x$ is isomorphic to some $\R^A$ consider any unconditional symmetric basis $\{e_\a\}_{\a\in A}$ of $\R^A$ and set $X_\a(y)=\Psi(e_\a,y)$ for any $y\in V$. We set ${\cal X}_x=\{X_\a,\;\a\in A\}$ and after a choice of such a set ${\cal X}_x$ for any $x\in M$, let be ${\cal X}=\dis\cup_{x\in M}{\cal X}_x$. By construction ${\cal X}$ satisfies {\rm{(Hi)}} and  { \rm{(Hiii)}} but  without (LB(s+2)). Given $x\in M$, with the previous notations,  we have  $||e_\a||_1=1$ and as $y\mapsto \Psi_y$ is a smooth field of continuous operators from $\R^A$ to $T_xM\equiv E$, we get  the property {\rm{(Hii)}} at $x$ after restriction of $V$ if necessary and also (LB(s+2)) at $x$ for {\rm{(Hiii)}}.\\

Now given  any  generating set of ${\cal D}$ which satisfies assumption (H), by application of part 1, we get the result.

\end{proof}

\subsection{Attainable set in infinite dimensional control theory for a family of vector fields}

Let $\cal X$ be a family of local vector fields which satisfies condition \rm{(Hi)} and  \rm{(Hii)} on a Banach manifold $M$. In our context  a {\bf controlled trajectory } of the {\bf  controlled system} associated to $\cal X$ is a curve $\g:I\ap M$ which is the integral curve of some vector field
\begin{eqnarray}\label{contZ}
Z(x,t,u)=\dis\sum_{k=1} ^pu_k(t)Z_k(x)
\end{eqnarray}
 associated to a family $\zeta=\{Z_k\}_{k=1,\cdots,p}\subset {\cal X}$ which satisfied the assumptions of Theorem 1  and  where  $u=(u_k)$  is a family of bounded curves  of class $L^1$ on some interval of $\R$ (see Theorem 1). In these conditions, $u$ is called {\bf the control} associated to $\g$. An {\bf admissible trajectory} is a curve $\g:[a,b]\ap M$ such that there exists  a finite  partition $a=t_0\leq \cdots\leq t_n$ such that   $\g: [t_i,t_{i+1}] \ap M $ is a controlled trajectory of the controlled system associated to $\cal X$   for $i=0,\cdots n-1$.

This context can be found in many papers (see for example: \cite{CH}, \cite{GXB}, \cite{XLG}, \cite{BH}, \cite{BP}, \cite{BBP},  \cite{Ro}). On the other hand, it is easy to see that any   $\cal X$-smooth piecewise curve is an admissible trajectory (see subsection 2.1). \\

According  to the classic context in control theory for a family $\cal X$  of vector  fields on $M$, the {\bf exact  attainable  set} ${\cal A}(x)$ of a point $x\in M$ is the set of points $y$ such  there exists an admissible trajectory $\g:[0,T]\ap M$ such that each $\g(0)=x$ and $\g(T)=y$. On the other end, the {\bf approximate  attainable  set} of $ x\in M$ is the closure $\bar{\cal A}(x)$ in $M$. \\

 \begin{rem}\label{orbatt}${}$\\ According to Proposition \ref{l1Xorbit}, if $\hat{\cal D}$ is integrable,  for any $\zeta$ as in (\ref{contZ}), on each maximal integral manifold $N$ which meets $V=\dis\cap_{k=1,\cdots,p}$Dom$(Z_k)$,   there exist vector fields $\tilde{Z}_k$, such that $(i_N)_*\tilde{Z}_k=Z_k$. So if we set $$\tilde{Z}(x,t,u)=\dis\sum_{k=1}^pu_k(t)\tilde{Z}_k(x)$$ then we have $(i_N)_* \tilde{Z}(t,u,.)=Z(t,u,.)$ and then we obtain that  each controlled trajectory with origin in $L$ is contained in $L$. In this case, if
  ${\cal O}(x)$ is the $\cal X$-orbit of $x$, we have the inclusions:
  $${\cal A}(x)\subset {\cal O}(x)\subset \bar{\cal A}(x)$$
 \end{rem}

In finite dimension,  we have ${\cal A}(x)={\cal O}(x)$ and it is well known  (from \cite {Su})  that ${\cal A}(x)$ is an {\bf immersed submanifold} of $M$ for any $x\in M$.


In our context, a corresponding result is given by the following Theorem:

 \begin{theor}\label{V}${}$\\
 Assume  that  the set $\hat{\cal X}$ (resp. the characteristic distribution  $\hat{\cal D}=l^{1}(\hat{\cal X})$) satisfies  the conditions (H) (see subsection \ref{loccond})(resp. (H') (see subsection \ref {stronginvol})) at any point $x\in M$.  Then $\hat{\cal D}$ is integrable.
 The exact  attainable  set ${\cal A}(x)$ of  any $x\in M$ is  dense in the maximal integral manifold $L(x)$ of $\hat{\cal D}$ through $x$ and the approximate attainable set $\bar{\cal A}(x)$ is  the closure of $L(x)$ in $M$  and also the closure of the $\cal X$-orbit of $x$. Moreover  if the distribution $\hat{\cal D}$ is closed $\hat{\cal A}(x)$  is a weak submanifold of $M$ for any $x\in M$.
 \end{theor}

The reader will find an illustration of this theorem in \cite{Ro} or in \cite{PS} (see also Example \ref{serpent}). Note that, if $\hat{\cal D}$ is finite dimensional, from Remark \ref{dimfin}, the assumptions of Theorem  \ref{V} are always satisfied and the distribution $\hat{\cal D}$ is closed. In this case the attainable set is exactly a $\cal X$-orbit. So in particular, when $M$ is finite dimensional we obtain Sussmann's result.\\

In finite dimension, to the  distribution $\cal D$ we can associate a chain of distributions
\begin{eqnarray}\label{chain}
{\cal D}^1={\cal D}\subset\cdots\subset {\cal D}^k\subset \cdots
\end{eqnarray}

\noindent where , for $k\geq 2$ $, {\cal D}^k$ is generated by the set ${\cal X}^k$ of local vector fields of type $[X_1,[\cdots [X_{k-1},X_k]\cdots ]$ where $X_1,\cdots, X_k$ belongs to $\cal X$.
The famous Theorem of Chow-Rashevsky asserts that if , for any $x\in M$, there exists $k$ such that
 ${\cal D}_x^k=T_xM$ then $M$ is the attainable  set of any point $x\in M$.\\

 Classically,    $\cal X$ is called  {\bf approximatively controllable} (resp. {\bf exactly controllable})  if  $\bar{\cal A}(x)=M$ (resp. ${\cal A}(x)=M$) for any $x\in M$. In order to to give an analogue of Theorem of Chow-Rashevsky  we have already associated to $\cal D$, a chain of distributions as in (\ref{chain}) (see subsection \ref{stronginvol}). As we have seen, if ${\cal X}$ satisfies condition (Hii) for some $s\in \N$, then the set ${\cal X}^k$ also satisfies (Hii)   for $s'=s-k$ and then, the $l^1$-characteristic distribution ${\cal D}^k=l^1({\cal X}^k)$ is well defined.  So we have the following version of Theorem of Chow-Rashevsky :

 \begin{theor}\label{VI}${}$\\
Assume  that  the set $\hat{\cal X}$ (resp. the characteristic distribution  $\hat{\cal D}=l^{1}(\hat{\cal X})$ satisfies  the conditions (H) (see subsection \ref{loccond})(resp. (H') (see subsection \ref {stronginvol}) at any point $x\in M$. Moreover, we suppose that   for any $x\in M$, there exists $k$ such that ${\cal D}^k_x$ is defined,
 and is dense in $T_xM$ (resp. ${\cal D}^k_x=T_xM$). Then $M$ is approximatively controllable (resp. exactly controllable).
 \end{theor}

In the previous Theorem, note that, according to the assumption we can have  controllability only if  the Banach manifold   $M$ is modeled on some $l^1(A)$ where $A$ is  a countable or uncountable set.  \\

We  say that  a distribution $\cal D$ on $M$ is finite co-dimensional  if for each $x$, the normed space ${\cal  D}_x$ is finite co-dimensional in $T_xM$. In this case ${\cal D}_x$ must be closed. In particular, finite co-dimensional $l^1$ distribution on  $M$ again imposes that $M$ is modeled on  $l^1(A)$ where $A$ is  a countable or uncountable set. In this case we have:

 \begin{cor}\label{VII} ${}$\\
Let $M$ be a Banach manifold modeled  on some $l^1(B)$. Consider  any set of vector fields $\cal X$ on $M$, which satisfies the conditions (H). If the characteristic distribution $\cal D$ is finite co-dimensional, then $M$ is foliated  by weak  Banach submanifolds of $M$ and each leaf is an $\cal X$-orbit.  Moreover, each attainable set is dense in such a leaf.

\end{cor}

\noindent \begin{proof}\so{\it Proof of Theorem \ref{V}}\\
 By application of  Theorem \ref{III} or Theorem \ref{III"}, we get the integrability of $\hat{\cal D}$. On one hand, for any $x\in M$, if $y$ belongs to ${\cal A}(x)$ as in (\ref{contZ}), the set $\zeta$ is finite,  it follows from Proposition \ref{varXS} that each integral curve of such a $Z$ is tangent to the leaf $L$ through $x$. On the other hand, from Proposition \ref{l1Xorbit},  if $y$ belongs to $L$, then $y$ is adherent to ${\cal A}(x)$.  According to Remark \ref{orbatt}, we have
 $${\cal A}(x)\subset L(x)\subset {\cal O}(x)\subset \bar{\cal A}(x)$$
 So, we get $\bar{L}(x)=\bar{\cal O}(x)=\bar{\cal A}(x)$.

 \noindent The last part is also a consequence of Theorem \ref{III} or Theorem \ref{III"}.

\end{proof}

\noindent\begin{proof}\so{\it Proof of Theorem \ref{VI}}${}$\\
From Theorem \ref{V} we know  that $\hat{\cal D}$ is integrable and, as Banach space is isomorphic to some $\R^A$. So by the same arguments as the ones used in the proof of Theorem \ref{III"} part 2, we have ${\cal D}^k_x\subset\hat{\cal D}_x$. It follows that $\hat{\cal D}_x$ is dense in $T_xM$ or equal to $T_xM$. The result is then a consequence of properties (ii) or (iii)  respectively of Theorem \ref{III'} or Theorem \ref{III"}.

\end{proof}

\noindent\begin{proof}\so{\it Proof of Corollary \ref{VII}}${}$\\
It is sufficient to prove that $\hat{\cal X}$ satisfies the condition (H) at each point $x\in M$. Given any $x\in M$, from our assumption we know that $\cal X$ satisfies the condition (H) at $x$. Take an unconditional symmetric basis $\{X_\a(x)\}_{\a\in A}$ such that $\{X_\a\}_{\a\in A}\subset {\cal X}_x$ and satisfies the  condition (LB(s+2)) for $s>0$. As  $\hat{\cal X}_x$ contains ${\cal X}_x$ and as ${\cal D}_x$ is finite co-dimensional, we can choose in $\hat{\cal D}_x$ a finite number $Y_1,\cdots Y_p$ such that $\{X_\a(x)\}_{\a\in A}\cup \{Y_1(x),\cdots Y_p(x)\}$ is an unconditional symmetric basis and $\{X_\a\}_{\a\in A}\cup \{Y_1,\cdots Y_p\}$ satisfies the condition  (LB(s+2)) for $s>0$. We then apply Theorem \ref{V}.
The last part can be shown as in the finite dimensional case (see \cite{Su})
 \end{proof}

\section{Proof of Theorem \ref{II}}\label{appendice}

In this last section,  we will use Theorem \ref{I} to give a proof of Theorem \ref{II}. \\
Recall that  $\xi = \{ X_\a, \; \a \in A\}$ is a family of vector fields defined on an open neighborhood $V$ of $x_0\in E$ and satisfies the condition (LB(s+2)) at $x_0$ and with the relation (\ref{LBs}) true for all $x\in V$.

\subsection {Maps $ \Gamma^{\tau}$  and
$\hat\Gamma^{\tau}$} \label{gamma1}

In this subsection we fix $\tau = (\t_\a)_{\a \in A}\in \R^A$. Let $B$ be any countable subset of $A$ which contains all the indexes $\a\in A$ such that $\t_\a\not=0$.  The set $B$ can be written as a sequence $\{\b_i,\; i\in \N\}\subset A$.   For the sake of simplicity, we then denote by $\t_i$ instead of $\t_{\b_i}$ the corresponding term of  $(\t_\a)_{\a \in A}$. With these notations we define
the sequence  $(\Gamma^{\tau}_i)_{i \in B}$  in the following way

\begin{itemize}
\item[\huge{.}] for $i = 1$
\begin{itemize}
\item[\normalsize{.}]  if  $ \t_1 = 0$ then $\Gamma^{\tau}_1(s) = 0$

\item[\normalsize{.}]if $ \t_1 \neq 0$ then $\Gamma^{\tau}_1(s)= \left \{
\begin{array}{lll}
 \displaystyle \frac{\t_1}{|\t_1|} & \textrm{if $s \in [0,|\t_1|[$} \\
&\\
0 & \textrm{othewise}\\
\end{array} \right.$
\end{itemize}
\item[\huge{.}]  for $i > 1$
\begin{itemize}
\item[\normalsize{.}] If $ \t_i = 0$ then $\Gamma^{\tau}_i(s) = 0$\\
\item[\normalsize{.}] If $ \t_i \neq 0$ then $\Gamma^{\tau}_i(s) = \left \{
\begin{array}{lll}
 \displaystyle \frac{\t_i}{|\t_i|} & \textrm{ if  $ s \in \displaystyle [\sum_{j=1}^{i-1}|\t_j|,\sum_{j=1}^{i}|\t_j|[$} \\
&\\
0 & \textrm{otherwise}
\end{array} \right.$
\end{itemize}
\end{itemize}

Now, for all $\a\not \in B$ we set $\G^\t_\a(s)=0$ for all $s\in \R$.\\

\noindent Finally, we define the families $\G^\t(s) $ and $\hat{\G}^\t(s)$ in the following way  by:

$\G^\t(s)=(\G^\t_\a(s))_{\a\in A}$ and  $\hat{\G}^\t=(\hat{\Gamma}^{\tau}_\a(s))_{\a\in A} =(\Gamma^{\tau}_\a(\| \tau \|_1-s))_{\a\in A}$

\noindent From this construction, it follows that:
$$\forall s \in
\R,\; (\Gamma^{\tau}_\a(s))_{\a \in A} \in
\R^A \textrm{ and } (\hat{\Gamma}^{\tau}_\a(s))_{\a \in A} \in
\R^A$$

\noindent  Now we consider the maps $\Gamma^{\tau}$ and
$\hat\Gamma^{\tau}$ defined in the following way:
\begin{displaymath}
\begin{tabular}{r r c l }
$\Gamma^{\tau} :$&$ \R$&$ \longrightarrow $&$\R^A$\\
&$s$&$\longmapsto$&$\Gamma^{\tau}(s) = (\Gamma^{\tau}_\a(s))_{\a\in
A}$
\end{tabular}
\end{displaymath}
\begin{displaymath}
\begin{tabular}{r r c l }
$\hat\Gamma^{\tau} :$&$ \R$&$ \longrightarrow $&$\R^A $\\
&$s$&$\longmapsto$&$\hat\Gamma^{\tau}(s) =
(\hat\Gamma^{\tau}_\a(s))_{\a\in A}$
\end{tabular}
\end{displaymath}

\begin{lem}\label{gamma2}${}$\\
 $\Gamma^{\tau}$ et $\hat\Gamma^{\tau}$ belongs to
 $L^1_b(\R)$.
\end{lem}

\subsection{Proof of the first part of Theorem \ref{II}}
In this section $x\in V$  and $\t=(\t_\a)_{\a\in A}\in \R^A$ are fixed. We consider any element  $\s=(\s_\a)_{\a\in A}$  of $\R^A$. We choose a countable subset $B$ of $A$ such that $B$ contains all the indexes $\a\in A$ such that  $t_\a\not=0$ and also all indexes $\b\in A$ such that $\s_\b\not=0$. Again the ordered set $B$ can be written as $B=\{\b_i,\;i\in \N\}$ and we then denote by $(\t_i)_{\b_i\in B}$ (resp. $(\s_i)_{\b_i \in B}$) the corresponding subsequence or $\t$ (resp. $\s$) and also we denote simply by $X_i$ the vector field $X_{\b_i}$ of $\xi$ for  all $\b_i\in B$.\\

\noindent With these notations,  for any $n\in \N$ and any $\s\in \R^A$, we set $\s^n=(\s_1,\cdots \s_n)\in \R^n$ and $\R^n$ is then considered as a subset of $\R^B\subset \R^A$

\begin{eqnarray}\label{psin}
 \psi^x_n(\tau^n) =\phi^{X_n}_{\t_n}\circ\ldots\circ\phi^{X_2}_{\t_2}\circ\phi^{X_1}_{\t_1}(x) \\ \nonumber
\hat{\psi}^x_n(\tau^n) =\phi^{X_1}_{\t_1}\circ\ldots\circ\phi^{X_n}_{\t_n}(x)
\end{eqnarray}

\begin{lem}\label{psin}${}$\\
With the previous notations, for each $n\in \N$, the map  $\psi^x_n$ is differentiable on $\R^n \cap B(0,\dis\frac{r}{k})$, and its differential is given by:

\begin{displaymath}
\begin{tabular}{c}
${D\psi_n}^x_{(\tau^n)}(\sigma^n)= $\\
$ \displaystyle
{D\phi^{X_n}_{\t_n}}_{(\psi_{n-1}(\tau^{n-1}))}\circ\ldots\circ
{D\phi^{X_1}_{\t_1}}_{(x)}[\sum_{p=1}^{n}\sigma_p
{D\phi^{X_1}_{-\t_1}}_{(\psi_1(\tau^1))}\circ\ldots \circ
{D\phi^{X_p}_{-\t_p}}_{(\psi_p(\tau^p))}(X_p(\psi_p(\tau^p)))]$
\end{tabular}
\end{displaymath}

Moreover, we have:
$$ \|{D\psi^x_{n}}_{(\tau^{n})}\|\leq
ke^{2k\|\tau\|_1}$$
for all $x\in V$ and any $n\in B$
\end{lem}

The proof of this Lemma is an elementary calculus by induction. For more details see \cite{La} chapter 5.\\

\noindent Afterwards, we will simply note, for any fixed $x\in V$:
$${D\psi_n}_{(\tau^n)}(\sigma^n)={D\phi^{X_n}_{t_n}}\circ\ldots\circ
{D\phi^{X_1}_{\t_1}}[\sum_{p=1}^{n}\sigma_p
{D\phi^{X_1}_{-\t_1}}\circ\ldots\circ
{D\phi^{X_p}_{-\t_p}}(X_p(\psi_p(\tau^p)))].$$

\noindent  We now define the following map:
\begin{eqnarray}\label{chin}
 \D\psi^x_n(\tau^n)=D\phi^{X_n}_{\t_n}\circ\ldots\circ D\phi^{X_1}_{\t_1}(x)\\ \nonumber
 \hat{ \D\psi}^x_n(\tau^n) =D\phi^{X_1}_{\t_1}\circ\ldots\circ D\phi^{X_n}_{\t_n}(x).
\end{eqnarray}

\noindent For these maps in the same way,  we obtain

\begin{lem}\label{chin}${}$\\
 For any fixed $x\in V$,  for each  $n\in \N$ the maps $\D\psi^x_n$ et $\hat{ \D\psi}^x_n$  are differentiable on $\R^n\cap B(0,\dis\frac{r}{k})$.
\end{lem}

\noindent \so{\it Now we are in situation to prove  part 1 of Theorem \ref{II}}\\

Let   $x_0\in V$ and $r>0$ be such that $
B_f(x_0,2r)\subset V$  and fix  $\tau = (\t_\a)_{\a \in A}\in \R^A$ such that  $\t\in B(0,\frac{r}{k}) \subset \R^A$. We fix  some countable subset $B\subset A$ which contains the set of indexes $\a$ such that $\t_\a\not=0$. As before the ordered set $B$ can be written $B=\{\b_i,\; i\in \N\}$  and  each $\t_{\b_i}$ with $\b_i \in B$  will be denoted $\t_i$. With these notations, we set
$$T =
\displaystyle \sum_{i \in \N} |\t_i| =\dis\sum_{\a\in A}|\t_\a|= \| \tau \|_1.$$
\textbullet$\quad$ Now we use  Theorem\ref{I}  with the following adaptations:
 $I=\R$ ,  $ \quad u =
\Gamma^{\tau}$,$\quad t_0 = 0$,
$\quad \delta$  a real number large enough
and $\quad T_0 = T$.

\noindent From Lemma \ref{gamma2} we have
$ \Gamma^{\tau} \in L^1_b(\R), \quad \textrm{ with }\quad
\|\Gamma^{\tau}\|_{\infty}=1$. As
$T < \displaystyle
\frac{r}{k}$, if we set  $I_0 = [-T,T]$ and  $U_0=B_0 = B(x_0,
r-kT)$, we get a flow $\Phi_{\Gamma^{\tau}}$, defined on  $I_0
\times U_0$. In particular, from Theorem\ref{I}
 the map
$\phi^{\xi}_{\tau}=\Phi_{\Gamma^{\tau}}(T,\;)$
is a $C^s$ diffeomorphism, and moreover, by construction, we get
$$\phi^\xi_\tau(x) = \displaystyle \lim_{n \to \infty}
\phi^{X_n}_{\t_n} \circ \ldots \circ
\phi^{X_1}_{\t_1}(x)=\displaystyle \lim_{n \to \infty}
\psi^x_n(\tau^n)$$
The same argument can be used to obtain the result concerning $\hat{\psi}^\xi_\t$.\\

\noindent \textbullet $\quad$ Now we prove that the inverse map of  $\phi^\xi_\tau$ is $\hat{\phi}^\xi_\tau$.\\
\begin{displaymath}
\begin{tabular}{ r c  l }
$\|\phi^\xi_\tau(\hat{\phi}^\xi_\tau(x))-x\|$ & $\le$ &
$\|\phi^\xi_\tau(\hat{\phi}^\xi_\tau(x)) -
\psi_n^{\hat{\phi}^\xi_\tau(x)}(\tau^n)\|+\|\psi_n^{\hat{\phi}^\xi_\tau(x)}(\tau^n) - x\|$\\
& $\le$ & $\|\phi^\xi_\tau(\hat{\phi}^\xi_\tau(x)) -
\psi_n^{\hat{\phi}^\xi_\tau(x)}(\tau^n)\|+\|\psi_n^{\hat{\phi}^\xi_\tau(x)}(\tau^n)-
\psi_n^{\hat{\psi}^x_n(-\tau^n)}(\tau^n)\|$\\
\end{tabular}
\end{displaymath}
At first, we have
$$\lim_{n \to \infty} \|\phi^\xi_\tau(\hat{\phi}^\xi_\tau(x)) -
\psi_n^{\hat{\phi}^\xi_\tau(x)}(\tau^n)\|= 0$$
So, it remains to show that the second term in the previous  majoration  converges to
 $0$ when  $n\ap \infty$.\\
The map  $x \longmapsto \psi_n^{x}(\tau^n)$ is of class
$C^1$ and its differential at $x$ is noting but
$\D\Psi_n(\tau^n)$.
So we have
$$\|\D\Psi_n^{x}(\tau^n)\| \le e^{kT}$$
So we obtain
 $$\|\psi_n^{\hat{\phi}^\xi_\tau(x)}(\tau^n) -
\psi_n^{\hat{\psi}^x_n(-\tau^n)}(\tau^n)\| \le
e^{kT}(\hat{\phi}^\xi_\tau(x)-\hat{\psi}^x_n(-\tau^n))$$
Finally, we get
$$\lim_{n \to \infty}
\|\psi_n^{\hat{\phi}^\xi_\tau(x)}(\tau^n) -
\psi_n^{\hat{\psi}^x_n(-\tau^n)}(\tau^n)\| = 0$$ which ends the proof of part $1$ of Theorem \ref{I}

\subsection{Proof of the second part of Theorem\ref{II}} \label{part2}

For any fixed $x \in U_0$,  we introduce the following notations:
\begin{eqnarray}\label{Psichi}
\psi^x_B(\tau)=\displaystyle \lim_{n \to \infty} \psi^x_n(\tau^n)=\phi^\xi_\tau(x)\\ \nonumber
\D\psi^x_B(\tau)=\displaystyle \lim_{n \to \infty}\D\psi^x_n(\tau^n)=D_2\Phi_{\Gamma^{\tau}}(T,x)
\end{eqnarray}

As a consequence  of Lemma \ref{psin}  and Lemma \ref{chin} we get:
\begin{lem}\label{continuitePsi}${}$\\
 $\psi^x_B$  and  $\D\psi^x_B$ are continuous maps on $\R^B\cap B(0,\dis\frac{r}{k})$
\end{lem}

 For each $\a\not\in B$ we can remark that $t_\a=0$ and so $\phi^{X_\a}_{t_\a}=Id$ and $D\phi^{X_\a}_{t_\a}=Id$ .
 So the previous limits (\ref{Psichi}) can be seen as an uncountable composition of maps of type $(\phi_\a)_{ \a\in A}$, evaluated at $x$, where only a countable subset of them are not equal to the identity.

\begin{nota}\label{uncountablecompo}${}$\\
 Given any $\a\in A$  we set   $\t^\a=(t_{\a'})\textrm{ with }\a'\in A, \;\a' \leq \a$. \\On the other hand for any $\a\in A$ we  consider the set
 $$B_\a=\{\b_i \textrm{ such that } \b_i\leq \a\}$$

 Considering the family of local diffeomorphisms associated to the family $\xi$ of vector fields we denote by:\\
 $\psi_\a^x(\t^\a)= \left \{
\begin{array}{lll}
\psi^x_n(\t^n) & \textrm{ if  } B_\a =\{\b_1,\cdots,\b_n\} \\
&\\
\psi^x_B &\textrm{ if  }  B_\a =B\\
\end{array} \right.$\\
${}$\\
  $\D\psi_\a^x(\t^\a)= \left \{
\begin{array}{lll}
\D\psi^x_n(\t^n) & \textrm{ if  }B_\a =\{\b_1\cdots\b_n\} \\
&\\
\D\psi^x_B &\textrm{if  }B_\a =B\\
\end{array} \right.$\\
${}$\\
 $\D\hat{\psi}_\a^x(\t^\a)= \left \{
\begin{array}{lll}
\D\hat{\psi}^x_n(\t^n) & \textrm{ if  }B_\a =\{\b_1\cdots\b_n\} \\
&\\
\D\hat{\psi}^x_B &\textrm{if  }B_\a =B\\
\end{array} \right.$\\

${}$\\
 $\psi^x(\t)=\psi_B^x(\t)=\phi^\xi_\t  \textrm{ and } \D\psi^x(\t)=\D\psi^x_B(\t)=D_2\Phi_{\Gamma^{\tau}}(T,x)=D\psi_\t^\xi(x).$\\
 \end{nota}

Given any $\s=(\s_\a)_{\a\in A}\in \R^A$, by taking for $B$ any countable set which contains the (countable) sets $\{\a\textrm{ such that }\t_\a\not=0\}$ and $\{\a\textrm{ such that }\s_\a\not=0\}$, from Lemma \ref{continuitePsi}   and Notations \ref{uncountablecompo} we get

\begin{lem}\label{continuPsia}${}$\\
The  map $\psi^x$ and $\D\psi^x$ are continuous on $B(0,\dis\frac{r}{k})$.
\end{lem}

\noindent \so{\it Now we can prove part {\it 2} of Theorem \ref{II}}\\

\noindent  We begin by proving that  $\psi^x$ is a $C^1$ map. We will use the following result of  \cite{D}, page 426:

\begin{prop}\label{existencediff}${}$\\
Let  $X$ et $Y$ be two  Banach spaces, $U\subset X$ an open set
and  $D $ a dense  vector subspace of  $X$. Consider  a continuous  map $f :\ U \to Y$
such that,
for all $(x,v)\in U\times X$, the  derivative   $f$ at $x$ in the direction
 $v$ denoted by $\partial_v{f}(x)$ exists. Moreover, assume
that there exists a continuous map $L : U \to L(X,Y)$ such that, for any $ (x,v)\in  U \bigcap D \times D$, we have
$\partial_v{f}(x)= L(x)(v)$. Then  $f \in {\cal{C}}^1(U,Y) \quad
\textrm {and } Df = L$
\end{prop}
We apply this result to the sets:

  $X = \R^A$, $\quad U
=B(0,\displaystyle \frac{r}{k})$  and $\quad Y = E$;

$ D = span\{e_\a, \a \in A\}$ where $e_\a = (\d^\a_\b)_{\b\in A}$,  where $\d^\a_\b=1$ if $\a=\b$ and $\d^\a_\b=0$ for $\a\not=\b$\\
 \indent (in fact, $ \{e_\a,\, \a\in A\}$ is the canonical basis of $\R^A$)

 the map $f$ is the map  $\psi^x$ on $B(0,\dis\frac{r}{k})\subset \R^A$

  $L$ is defined in the following way:

\noindent for $\sigma = (\sigma_\a)_{\a \in A}\in
B(0,\displaystyle \frac{r}{k}) \subset \R^A$ :
$L(\tau)(\sigma)=\D\psi^x(\tau)
[\displaystyle\sum_{\a\in A}\sigma_\a
{\D\hat{\psi}^x( -\t^\a)}[X_\a(\psi_\a(\tau^\a)))]$.

$\bullet$ It is clear that $D$ is a dense set in $\R^A$.

$\bullet$ The continuity of $\psi^x$ follows from Lemma \ref{continuPsia}.

$\bullet$ Now we prove that $ \quad  \forall \tau  \in
\displaystyle B(0, \frac{r}{k}) \cap D,\forall \sigma \in D,\quad
\partial_{\tau}{\psi^x}(\sigma)= A(\tau)(\sigma)$\\

Let be $(\tau ,\sigma)\in \displaystyle B(0,
\frac{r}{k}) \cap D \times D .$ So we have
$$\tau = (\t_{\a_i})_{i=1,\cdots,p}, \quad \textrm{with} \quad \|\tau\|_1 < \frac{r}{k}$$
$$\sigma = (\s_{\b_j})_{j=1\cdots,q}\quad \textrm{with} \quad ||\s||_1<\frac{r}{k}$$
\indent The family $\{(e_{\a_i})_{i=1,\cdots,p},\; (e_{\b_j})_{j=1\cdots,q}\}$ can be put in an ordered family $\{e_{\a_l}\}_{l=1\cdots n}$ with $n\leq \inf (p,q)$.
So we can consider that $\t$ and $ \s$ belong to $span\{e_{\a_1},\cdots,e_{\a_n}\}$. For simplicity we denote by $\t_i$ (resp. $\s_i$) the component  of $\t$, (resp. $\s$) on $e_{\a_i}$ and $X_i$ instead of $X_{\a_i}$, for $i=1,\cdots, n$.
Now, for any $\l\in \R$ and $\l\not=0$, we have:
\begin{displaymath}
\begin{tabular}{r  c  l}
$\psi^x(\tau +\lambda\sigma) - \psi^x(\tau )$ & $ = $ & $\psi^x_n(\tau^n
+\lambda\sigma^n)-\psi_n(\tau^n )$\\
&&\\
& $ = $ & $\displaystyle\lambda
{D\phi^{X_n}_{\t_n}}\circ\ldots\circ
{D\phi^{X_1}_{\t_1}}[\sum_{p=1}^{n}\sigma_p
{D\phi^{X_1}_{-\t_1}}\circ\ldots$\\
&&\\
&  & $\hfill \quad \quad \quad \ldots\circ
{D\phi^{X_p}_{-\t_p}}(X_p(\psi_p(\tau^p)))] + o(\lambda\sigma^n)$
\end{tabular}
\end{displaymath}
so:
$$\displaystyle \partial_{\tau}{\Psi}(\sigma)= \lim_ { \lambda
\rightarrow 0} \frac {\Psi(\tau +\lambda\sigma) - \Psi(\tau
)}{\lambda}= L(\tau)(\sigma)$$

$\bullet$  the continuity of $\t\ap L(\t)$ :

 \noindent Now we consider the following map:\\
\indent \begin{tabular}{r c c l }
${\cal R}:$&$B(0, \r) $&$\longrightarrow$&${\cal{L}}(\hat{\cal D}_x)$\\
&$\t$&$\longmapsto$&${\cal R}(\t)$\\
\end{tabular}

defined by
 ${\cal R}(\t)(\dis\sum_{\a\in A}\s_\a X_\a(x))=\dis\sum_{\a\in A}\s_\a\D\hat{\psi}^x_\a((-\tau)^\a)[X_i(\psi^x_\a(\tau^\a))]$

 Note that  from Lemma \ref{psin}, we have
 $$||{\cal R}(\t)(\dis\sum_{\a\in A}\s_\a X_\a(x))||\leq ke^{k||\t||_1}||\s||_1$$
 So ${\cal R}(\t)$ is a continuous linear map.
 On the other hand, we can write

 \begin{eqnarray}\label{decompositiondif}
 L(\t)(\s)=\D\psi^x(\tau)\circ{\cal R}(\t)(\dis\sum_{\a\in A}\s_\a X_\a(x))
\end{eqnarray}

The proof of the following Lemma can be found in \cite{La} chapter 5:
\begin{lem}\label{contR}${}$\\
The map $\t \mapsto {\cal R}(\t)$ is continuous on $B(0,\frac{r}{k})$.
\end{lem}
From this lemma and Lemma \ref{continuPsia}, it follows that $\t\mapsto L(\t)$ is continuous.\\

So we obtain that $\psi$ is $C^1$on $\displaystyle
B(0, \frac{r}{k})$.\\

To prove that $\psi$ is of class $C^s$ for $s\geq 2$, as classically we use the fact that:
$$(\Phi_{\G^\t}(t,x),D_2\Phi_{\G^\t}(T,x),\cdots,D^s_2\Phi_{\G^\t}(T,x))$$
is the flow of an adapted  vector field $\hat{Z}^s$ on an open set of the Banach space $E\times {\cal L}(E,E)\times \cdots\times {\cal L}^s(E,E)$ and proceed by induction.


\end{document}